%% file: source.tex
\documentclass{elsarticle}
\input{preamble}

\begin{document}

\begin{frontmatter}



\title{Worst-case multi-objective error estimation and adaptivity}


\author[TUE1,TUE2]{E.H.~van~Brummelen}
\author[IBM1]{S.~Zhuk}
\author[TUE1]{G.J.~van~Zwieten}

\address[TUE1]{%
  Eindhoven University of Technology, Department of Mechanical Engineering,
  P.O. Box 513, 5600 MB Eindhoven, The Netherlands}
\address[TUE2]{%
  Eindhoven University of Technology, Department of Mathematics \& Computer Science,
  P.O. Box 513, 5600 MB Eindhoven, The Netherlands}
\address[IBM1]{%
  IBM Research, Server 3, IBM Tech. Campus, Damastown, Dublin 15, Ireland}

\begin{abstract}
This paper introduces a new computational methodology for determining
a-posteriori multi-objective error estimates for finite-element approximations, and for constructing corresponding (quasi-)optimal adaptive refinements of finite-element spaces. As opposed to the classical goal-oriented
approaches, which consider only a single objective functional, the presented methodology applies to general closed convex subsets of the dual space and constructs a worst-case error estimate of the finite-element approximation error.
This worst-case multi-objective error estimate conforms to a dual-weighted residual, in which the dual solution is associated with an approximate supporting functional of the objective set
at the approximation error. We regard both standard approximation errors and data-incompatibility errors associated with incompatibility of boundary data with the trace of the finite-element space. Numerical experiments are presented to demonstrate the efficacy of applying the proposed worst-case multi-objective error in adaptive refinement procedures.

\end{abstract}


\end{frontmatter}



\section{Introduction}
Goal-oriented {\em a-posteriori\/} error-estimation and adaptivity has emerged over the past years as an effective
methodology to solve measurement problems in computational science and engineering. As opposed to the classical
norm-oriented adaptive methods that originated from the seminal work of Babu\v{s}ka and Rheinboldt~\cite{Babuska:1978lq,Babuska:1978tg},
which aim to construct finite-element approximations that are optimal in an appropriate norm, goal-oriented adaptive methods
seek to construct finite-element approximations that yield optimal approximations of a particular functional of the
solution. Goal-oriented methods date back to the pioneering work of Becker and Rannacher~\cite{BeckerRannacher1996,BeckerRannacher2001},
Oden and Prudhomme~\cite{PrudhommeOden1999,Oden:2001ss} and Giles, S\"uli, Houston and Hartmann~\cite{GilesSuli2002,HoustonSuli2001,HartmannHouston2002a,HartmannHouston2002b}.
It is noteworthy, however, that the use of dual solutions in a-posteriori error estimates, which forms the basis of goal-oriented error estimation,
had already been pursued before by Johnson, Eriksson and Hansbo~\cite{Eriksson:1991dz,Johnson:1990tk,Johnson:1992fk}.
Goal-oriented adaptivity has been successfully applied to a wide variety
of problems including compressible and incompressible flow problems~\cite{HartmannHouston2002b,PrudhommeOden2003,Rannacher2009},
elasticity and plasticity~\cite{Stein:2004kx,Rannacher:1997hl,Rannacher:1999to}, variational inequalities~\cite{Suttmeier:2001gd},
fluid-structure interaction~\cite{Zee:2011uq,Fick:2010kx,Richter:2012uq,LarsonBengzon2010}, free-boundary problems~\cite{Zee:2010fl,Zee:2010uq},
phase-field models~\cite{Zee:2011vn,Wick:2015bx}, and kinetic equations~\cite{Hoitinga2011AIP}. This overview is in fact far from exhaustive, and many
more applications and contributions can be found in the literature. It is notable, however, that despite the success of goal-oriented methods in
many applications, convergence and optimality theory is still rudimentary; see~\cite{Mommer:2009yb,Feischl:2016vp}.

In addition to goal-oriented adaptivity of the mesh width~($h$) or local order~($p$) of finite-element approximation spaces,
goal-oriented approaches have also been considered in the context of model adaptivity. In these approaches, the error indicator is applied to
systematically decide between a simple coarse model and a complex sophisticated model, in such a manner that the sophisticated model is
only applied in regions of the domain that contribute most significantly to the objective functional under consideration. For examples
of goal-oriented model adaptivity, we refer to~\cite{OdenVemaganti2000} for application of goal-oriented model adaptivity to heterogeneous
materials, to~\cite{Bauman:2009kp} for goal-oriented atomistic/continuum adaptivity in solid materials, and to~\cite{Opstal:2015kr}
for goal-oriented adaptivity between a boundary-integral formulation of the Stokes equations and a PDE formulation of the Navier-Stokes equations.

The premise in goal-oriented strategies is that interest is restricted to a single functional output, rather than all details of the solution
of the system under consideration. Although there are indeed many applications in which interest is restricted to a single scalar output of a system,
there are in fact many more in which the object of interest is not so narrowly defined. There are for instance many applications in
which multiple scalar outputs are of interest, e.g. the aerodynamic drag, lift and torques that are exerted by a flow on an immersed object~\cite{Hartmann:2008db}.
Another example is furnished by the many cases in which one is concerned with the solution in a particular region of the problem domain. This scenario
also covers boundary-coupled multiphysics problems, e.g. fluid-structure interaction~\cite{Zee:2011uq},
in which two subsystems interact via a mutual interface and only one of the subsystems is of importance, and
free-boundary problems~\cite{Zee:2010fl,Zee:2010uq}, in which interest extends to the geometry of the domain only
and not to the solution of the underlying partial differential equation, or vice versa.
A similar situation arises for volumetrically-coupled multiphysics problems, e.g. electro-thermo-mechanical problems~\cite{LarsonBengzon2008},
if significance is assigned to one field only. Another pertinent class of problems that are incompatible with the premise of goal-oriented strategies,
pertains to applications in which one is interested in the
flux functional on a part of the boundary that is subjected to essential boundary conditions or, reciprocally, the trace of the
solution on a part of the boundary that is furnished with natural boundary conditions. In all of the aforementioned cases, interest does not
extend to all details of the solution, but it is also not restricted to a single quantity. Instead, the objective set in these
cases corresponds to a proper non-singleton subset of the dual space, and not to an individual element of the dual space as in
goal-oriented approaches.

The purpose of this paper is to introduce a new computational methodology for determining a-posteriori {\em multi-objective\/} error estimates for
finite-element approximations, and for constructing corresponding (quasi-)optimal adaptive refinements of finite-element spaces. Such
multi-objective error estimation and adaptivity has only received scant consideration so far; see~\cite{Estep:2005dp,Hartmann:2008db}.
As opposed to the aforementioned approaches, our methodology applies to generic possibly infinite-dimensional closed convex subsets of the dual space.
The presented methodology relies on the construction of a worst-case error estimate, viz. the supremum of the error over the objective set.
This worst-case error corresponds to the value of the supporting function of the objective set at the finite-element approximation error. We then construct
an approximate supporting functional of the objective set at the approximation error, and determine an approximate dual solution for this supporting
functional to form an error estimate in dual-weighted-residual form~\cite{BeckerRannacher1996}. The corresponding error estimate represents
a {\em worst-case multi-objective error estimate} for the finite-element approximation on the objective set. In this work, we regard both standard approximation
errors and data-incompatibility errors associated with incompatibility of boundary data with the trace of the finite-element space.


The presented worst-case multi-objective error estimate is based on a technique for estimating linear functionals of solutions of equations with unbounded closed linear operators set in Hilbert spaces~\cite{Zhuk2009d}. This technique applies to general linear operators with non-trivial null-spaces and possibly non-closed ranges, and it is based on an extension of Young--Fenchel duality to closed
unbounded linear operators~\cite{Zhuk2012AMO,Zhuk2009d}. It is noteworthy that this technique for estimating a linear functional of a solution of an operator equation, which we consider in this paper in the context of worst-case multi-objective error estimation, yields a natural generalization of
adaptive-state-estimation or so-called data-assimilation procedures (such as online sequential filters~\cite{Tirupathi:2016fk,ZhukSISC13,Tchrakian:2014rc} or offline variational approaches~\cite{Herlin:2012gp}) to a wide class of linear and nonlinear partial differential equations, enabling incorporation of a-posteriori knowledge (for instance, sensor readings) into the error estimate rendering it less conservative~\cite{ZhukSISC13,Zhuk2009c}.

The remainder of this paper is organized as follows. Section~\SEC{ModProb} presents a model problem that serves as a concrete setting
for the worst-case multi-objective error-estimation procedure. In this section we also review standard goal-oriented error estimation.
Section~\SEC{multi-object-error} introduces worst-case multi-objective error estimation in an abstract Hilbert-space setting.
We consider both standard finite-element-approximation errors and errors engendered by data incompatibility. In Section~\SEC{numer-exper},
we conduct numerical experiments to illustrate the properties of the worst-case multi-objective error-estimation procedure, and we consider
its application in an adaptive-refinement procedure. Section~\SEC{Concl} presents concluding remarks.

\section{Model problem and preliminaries}
\label{sec:ModProb}
To provide a concrete background for the presented worst-case multi-objective error-estimation
and adaptivity processes, and to elucidate aspects related to data incompatibility and
boundary conditions, we consider the standard Dirichlet--Poisson problem. A general
setting is presented in Section~\ref{sec:multi-object-error}.

\subsection{Problem statement}
\label{sec:ProbStat}
Consider a bounded open subset $\Omega\subset\IR^d$ ($d=2,3$)
with boundary $\partial\Omega$.
We assume throughout that~$\Omega$ is a Lipschitz domain, with
boundary~$\Gamma:=\partial\Omega$. We consider the standard Dirichlet--Poisson problem:
\begin{subequations}
\label{eq:BVP}
\begin{alignat}{2}
-\Delta{}u&=f&\quad&\text{in }\Omega
\label{eq:DiffEq}
\\
u&=g&\quad&\text{on }\Gamma
\label{eq:BCs}
\end{alignat}
\end{subequations}
where $f:\Omega\to\IR$ and $g:\Gamma\to\IR$ are exogenous functions.

To accommodate the weak formulation of~\EQ{BVP}, let $H^1(\Omega)$ denote the
Sobolev space of square-integrable functions with square-integrable (in Lebesgue sense) weak derivatives, equipped with the standard scalar product, norm and seminorm, here denoted
by $(\cdot,\cdot)_{H^1(\Omega)}$, $\|\cdot\|_{H^1(\Omega)}$ and~$|\cdot|_{H^1(\Omega)}$,
respectively. We denote
by~$\gamma:H^1(\Omega)\to{}L^2(\Gamma)$ the trace operator, defined via the continuous
extension of the operator that restricts continuous functions on~$\Omega$ to the
boundary~$\Gamma$. The trace space $H^{1/2}(\Gamma):=\gamma{}H^1(\Omega)$ corresponds to
a proper subspace of~$L^2(\Gamma)$. We denote by $H^1_0(\Omega)$ the kernel of~$\gamma$,
viz. the subclass of functions in~$H^1(\Omega)$ that vanish on~$\Gamma$ in the trace sense.
In addition,
we denote by $\smash[b]{\ell_{(\cdot)}:H^{1/2}(\Gamma)\to{}H^1(\Omega)}$ a suitable right-inverse of~$\gamma$.
Such a right inverse is generally called a lift operator. Lift operators
are non-unique. A particular example is the (harmonic) Moore--Penrose lift associated with the $H^1$\nobreakdash-semi norm:
\begin{equation}
\label{eq:harmonicl}
g\mapsto\argmin\big\{|v|_{H^1(\Omega)}:v\in{}H^1(\Omega),\gamma{}v=g\big\}
\end{equation}
The boundary value problem~\EQ{BVP} can now be condensed into the weak formulation: given $f\in{}L^2(\Omega)$ and $g\in{}H^{1/2}(\Gamma)$, find
\begin{equation}
\label{eq:weakform}
u\in\ell_g+H^1_0(\Omega):
\qquad
a(u,v)=b(v)\qquad\forall{}v\in{}H^1_0(\Omega)
\end{equation}
where the bilinear form $a:[H^1(\Omega)\times{}H^1(\Omega)]^{\star}$ and linear form $b\in[H^1(\Omega)]^{\star}$ are defined by:
\begin{equation}
\label{eq:a(u,v)b(v)}
a(u,v)=\int_{\Omega}\grad{}u\cdot\grad{}v,
\qquad
b(v)=\int_{\Omega}fv
\end{equation}
with $(\cdot)^{\star}$ the topological dual of~$(\cdot)$.
Weak problem~\EQ{weakform} with non-homogeneous boundary data~$g$ and  lift $\ell_g$ admits a
reinterpretation as:
\begin{equation}
\label{eq:weakform2}
u_0\in{}H^1_0(\Omega):\qquad
a(u_0,v)=b(v)-a(\ell_g,v)
\qquad\forall{}v\in{}H^1_0(\Omega)
\end{equation}
followed by the additive composition $u=\ell_g+u_0$.

\subsection{Galerkin finite-element approximations}
\label{sec:FEM}
Consider a strictly decreasing sequence of mesh parameters
$\mathcal{H}=(h_i)_{i\in\IZ_{\geq{}0}}$ with only accumulation point zero, and a
corresponding sequence of asymptotically dense
nested (finite-element) approximation
spaces
$V^{h_0}\subset{}V^{h_1}\subset{}V^{h_2}\subset\cdots\subset{}H^1(\Omega)$.
With each $h\in\mathcal{H}$ we associate a discrete lift operator $\smash[b]{\ell^h_{(\cdot)}:H^{1/2}(\Gamma)\to{}V^h}$
such that $\smash[b]{\gamma(\ell^h_g)}$ is close to $g$ in some appropriate sense.
Denoting by $V^h_0$ the intersection of $V^h$ and $H^1_0(\Omega)$, composed of approximation functions in~$V^h$ that vanish on the boundary,
we associate with each~$h\in\mathcal{H}$ the Galerkin approximation problem:
\begin{equation}
\label{eq:GalerkinPrimal}
u^h\in{}\ell^h_g+V^h_0:\qquad
a(u^h,v^h)=b(v^h)
\qquad\forall{}v^h\in{}V^h_0
\end{equation}
Defining the {\em residual\/}
operator~$r:H^1(\Omega)\to{}[H^1(\Omega)]^{\star}$ according to
\begin{equation}
\label{eq:residual}
\text{for all }(u,v)\in{}H^1(\Omega)\times{}H^1(\Omega):
\qquad
\langle{}r(u),v\rangle=b(v)-a(u,v)
\end{equation}
we note for later reference that the residual corresponding to the Galerkin
approximation~$u^h$ in~\EQ{GalerkinPrimal} is orthogonal with respect to
the test functions~$V^h_0$, i.e., $\langle{}r(u^h),v^h\rangle$ vanishes for all $v^h\in{}V^h_0$.


\subsection{Goal-oriented error estimation}
\label{sec:GOEE}
Consider a functional $\jmath\in{}[H^1(\Omega)]^{\star}$ and suppose that interest is restricted to $\jmath(u)$ with $u$ according to~\EQ{weakform}, rather than to the all details of the
solution to~\EQ{weakform}. A problem of this type is called a
{\em measurement problem} and $\jmath$ is generally referred as a {\em goal functional}, {\em target functional} or~{\em objective functional}.
By virtue of the Riesz representation theorem, the functional $\jmath(\cdot)$ can generally
be written in the form:
\begin{equation}
\label{eq:jform}
\jmath(\cdot)=a(\cdot,\vartheta)
\end{equation}
with $\vartheta\in{}H^1(\Omega)$. The trace of $\vartheta$ is associated with
a weak (so-called {\em variationally consistent\/}) formulation of the
flux~$\partial_nu$~\cite{Brummelen:2012fk,MelboKvamsdal2003}.
Indeed, if $u$ satisfies~\EQ{BVP}, then the
following integration by parts identity holds for all $\vartheta\in{}H^{1}(\partial\Omega)$:
\begin{equation}
\label{eq:fluxSeq1}
\int_{\partial\Omega}(\gamma\vartheta)\partial_nu
=
\int_{\Omega}\grad\vartheta\cdot\grad{}u
+
\int_{\Omega}\vartheta\Delta{}u
=
a(u,\vartheta)-b(\vartheta)
=:
\langle\gamma\vartheta,\partial_nu\rangle
\end{equation}
The functional $\langle\gamma\vartheta,\partial_n(\cdot)\rangle$ is affine and consists
of a linear part, $a(\cdot,\vartheta)$ and
a translation, $b(\vartheta)$.
The translation depends exclusively on data of the problem, and
is therefore in principle explicitly computable.
Only the linear part $a(u,\vartheta)$
of $\langle\gamma\vartheta,\partial_nu\rangle$ depends on the solution.

We associate with~\EQ{weakform} and the goal functional~$\jmath$
the following {\em dual problem}:
\begin{equation}
\label{eq:dual0}
z\in{}H^1_0(\Omega):
\qquad
a(v,z)=\jmath(v)
\qquad\forall{}v\in{H}^1_0(\Omega)
\end{equation}
In this context, problem~\EQ{weakform} is referred to as the {\em primal problem}.
It is to be noted that the function~$z$ corresponds to the Riesz representation
of $\jmath$ in~$H^1_0(\Omega)$ with inner product~$a(\cdot,\cdot)$. In view
of the representation of~$\jmath$ in~\EQ{jform}, the dual problem~\EQ{dual0} can alternatively be conceived of as a projection from $H^1(\Omega)$ into~$H^1_0(\Omega)$.
We will denote this projection by $\pi_0(\cdot)$. For any $\vartheta\in{}H^1(\Omega)$,
the complement $\vartheta-\pi_0\vartheta$ corresponds to the weakly harmonic lift of the trace
of~$\vartheta$. To corroborate this assertion, we note that by~\EQ{jform} and~\EQ{dual0} it holds that $a(v,\vartheta-\pi_0\vartheta)=0$ for all $v\in{}H^1_0(\Omega)$, which implies
that $\vartheta-\pi_0\vartheta$ is weakly harmonic.
Evidently, it moreover holds that~$\gamma(\vartheta-\pi_0\vartheta)=\gamma\vartheta$.

Consider the Galerkin approximation $u^h=u_0^h+\ell_g^h\in{}V^h$ to $u$
according to~\EQ{GalerkinPrimal}. The dual solution enables us to express
the error in the objective functional, $\smash[t]{\jmath(u)-\jmath(u^h)}$,
without direct reference to the error in the approximation.
The error representation follows from the following sequence of identities:
\begin{equation}
\label{eq:erf}
\begin{aligned}
\jmath(u)-\jmath(u^h)
&=
\jmath(\ell_g-\ell^h_g)+\jmath(u_0-u^h_0)
\\
&=\jmath(\ell_g-\ell^h_g)+a(u_0-u_0^h,z)
\\
&=\jmath(\ell_g-\ell^h_g)-a(\ell_g-\ell_g^h,z)+b(z)-a(u^h,z)
\\
&=a(\ell_g-\ell^h_g,\vartheta-\pi_0\vartheta)+\langle{}r(u^h),\pi_0{\vartheta}-v^h\rangle
\end{aligned}
\end{equation}
for all $v^h\in{}V^h_0$.
The first identity follows from the partitions $\smash[tb]{u=u_0+\ell_g}$
and~$\smash{u^h=u^h_0+\ell_g^h}$.
The second identity holds on account of $\smash[tb]{u_0-u_0^h\in{}H^1_0(\Omega)}$ and the dual
probem~\EQ{dual0}. The third identity follows from the primal problem in the form~\EQ{weakform2} and a partition of zero. The fourth identity follows from~\EQ{residual} and~\EQ{jform} and a rearrangement of terms, in combination with
the property~$\smash[t]{\langle{}r(u^h),v^h\rangle=0}$ for all $\smash[t]{v^h\in{}V^h_0}$ of the Galerkin approximation. It is to be noted that the
error-representation formula, viz., the ultimate expression in~\EQ{erf},
depends on the approximation~$u^h$ only via the corresponding residual.

The first term in the error representation can be conceived of as the error
induced by a discrepancy in the data of the weak formulation~\EQ{weakform}
and its Galerkin finite-element approximation~\EQ{GalerkinPrimal}. Typically, this discrepancy
results from an incompatibility of the boundary data, viz. $g\notin\gamma{}V^h$.
The second term represents the standard
{\em Dual-Weighted Residual (DWR)\/} contribution to the error in the objective
functional~\cite{BeckerRannacher1996}. In particular,
the error representation conveys that
the contribution of the boundary-data error, $\smash{\ell_g-\ell_g^h}$, to the
error in the objective functional vanishes if the representation function
of the objective functional, $\vartheta$, is compactly supported in~$\Omega$,
i.e. if $\smash{\vartheta\in{}H^1_0(\Omega)}$. It then holds
that the projection $\pi_0{\vartheta}$ coincides with~$\vartheta$. Therefore,
the DWR term separately provides an error representation if the boundary conditions
are compatible or if the trace of the representation function
of the objective functional vanishes.

\paragraph{Remark}
Two alternative error-representation formulas can be derived, in
addition to the final expression in~\EQ{erf}. To this end, let us
consider the Galerkin approximation of the dual problem~\EQ{dual0} in
the finite-element approximation space~$V^h$:
\begin{equation}
\label{eq:DualGalerkin}
z^h\in{}V^h_0:\qquad{}a(v^h,z^h)=\jmath(v^h)
\qquad\forall{}v^h\in{}V^h_0
\end{equation}
Denoting by $\rho:H^1(\Omega)\to{}[H^1(\Omega)]^{\star}$ the {\em dual residual\/},
\begin{equation}
\label{eq:dualresidual}
\text{for all }(u,v)\in{}H^1(\Omega)\times{}H^1(\Omega):
\qquad
\langle{}\rho(v),u\rangle=\jmath(u)-a(u,v)
\end{equation}
we infer that the residual of the Galerkin approximation~$z^h$ to the dual solution
according to~\EQ{DualGalerkin} satisfies the orthogonality relation~$\langle{}\rho(z^h),v^h\rangle=0$
for all $v^h\in{}V^h$. Therefore,
\begin{equation}
\label{eq:seqID2}
\begin{aligned}
\jmath(u)-\jmath(u^h)
&=
\jmath(\ell_g-\ell_g^h)+\jmath(u_0-u_0^h)
\\
&=
a(\ell_g-\ell_g^h,\vartheta)
+
\jmath(u_0)-a(u_0,z^h)+\big(a(u_0,z^h)-a(u^h_0,z^h)\big)
\\
&=
a(\ell_g-\ell_g^h,\vartheta-z^h)+
\langle{}\rho(z^h),u_0-v^h\rangle
\end{aligned}
\end{equation}
for all $v^h\in{}V^h_0$. Hence, the error in the goal functional, $\jmath(u)-\jmath(u^h)$,
can alternatively be expressed as the sum of a data-incompatibility contribution and
the duality pairing of the dual residual associated with
the Galerkin approximation of the dual solution, $\rho(z^h)$, and an interpolation error
in the homogeneous part of the primal solution in~$\smash[b]{V^h_0}$, i.e., $\smash[b]{u_0-v^h}$ for some $\smash[b]{v^h\in{}V^h_0}$. It is to be noted, however, that the first term in the ultimate
expression in~\EQ{seqID2} generally vanishes only if $\smash[t]{\ell_g-\ell_g^h=0}$ and not if $\vartheta\in{}H^1_0(\Omega)$, as opposed to the data-incompatibility
contribution in~\EQ{erf}.

In the absence of data-incompatibility errors,
we can also infer the following {\em symmetric error representation}
from~\EQ{erf} and~\EQ{seqID2}:
\begin{equation}
\label{eq:symerr}
\jmath(u)-\jmath(u^h)
=\tfrac{1}{2}\langle{}r(u^h),z-w^h\rangle+\tfrac{1}{2}\langle\rho(z^h),u-v^h\rangle
\end{equation}
for all $v^h,w^h\in{}V^h_0$.
For nonlinear problems, i.e.,
for semi-linear forms $a(\cdot,\cdot)$ and nonlinear goal functionals~$\jmath(\cdot)$,
the error representations generally deviate from the exact error on account of linearization errors.
The primal and dual error estimates exhibit quadratic remainders, i.e.,
the error estimates deviate from $\jmath(u)-\jmath(u^h)$ by $O(\|u-u^h\|_{H^1(\Omega)}^2)$
as $\|u-u^h\|_{H^1(\Omega)}\to{}0$. The symmetric error representation has the advantage
that it exhibits a third order remainder; see~\cite{BeckerRannacher2001}.

The worst-case multi-objective error-estimation process proposed in this paper and
the corresponding adaptive procedure are based on the error
representation~\EQ{erf}.\closeremark

\section{Worst-case multi-objective error estimation}
\label{sec:multi-object-error}
In many applications, interest does not extend to all details of a solution
to a boundary-value problem,
but it is also not restricted to a single objective functional.
Instead, one can envisage situations in which all functionals in a proper subclass
 of the dual space are of concern. Accordingly, to assess the accuracy
 of an approximation to the solution, an estimate of the associated error in
 all functionals in the considered subclass is required. In this section we propose
 a methodology for constructing such a {\em multi-objective error estimate\/},
 based on the worst-case error.

\subsection{Worst-case multi-objective error estimation without data incompatibility}
\label{sec:multi-object-error-compatible}
We formulate the worst-case multi-objective error estimate in an abstract setting encompassing
the model problem presented in Section~\SEC{ModProb}. For simplicity, we first restrict
ourselves to error estimates without data-incompatibility contributions (cf. Section~\SEC{GOEE}), so that the DWR error-representation formula for objective functionals applies directly. The extension to error estimates with data incompatibility is treated in Section~\ref{sec:multi-object-error-incompatible}.

Let $H$ denote a Hilbert space with inner product $(\cdot,\cdot)_H$ and induced
norm $\|\cdot\|_H$. The dual space~$H^{\star}$ of~$H$ consists of all continuous linear functionals on~$H$, equipped with the dual norm
$\|\cdot\|_{H^{\star}}=\sup\{\langle\,\cdot\,,v\rangle:v\in{}H,\|v\|_H=1\}$.
We consider a continuous and coercive bilinear form $a\in[H\times{}H]^{\star}$, a continuous linear
form $b\in{}H^{\star}$, and the problem:
\begin{equation}
\label{eq:P}
u{}\in{}H:\qquad{}a(u,v)=b(v)\qquad\forall{}v\in{}H
\end{equation}
The Lax-Milgram Lemma (see, for instance, \cite[Thm 2.7.7]{BrennerScott2002}) asserts that
under the aforementioned hypotheses on~$a$ and~$b$, problem~\EQ{P} is well posed and that the solution $u$ to~\EQ{P} satisfies the a-priori estimate $\|u\|_H\leq{}C\|b\|_{H^{\star}}$ for some positive constant $C>0$.
Denoting by $H^h\subset{}H$ a proper subspace of~$H$, the Galerkin approximation
of~\EQ{P} in~$H^h$ writes:
\begin{equation}
\label{eq:Ph}
u^h{}\in{}H^h:\qquad{}a(u^h,v^h)=b(v^h)\qquad\forall{}v^h\in{}H^h
\end{equation}

Assuming that interest is restricted to a closed and convex class of objective
functionals $\DD{}\subset{}H^{\star}$, it is appropriate to assess
the accuracy of the approximation~$u^h$ according to~\EQ{Ph} not
with reference to its deviation from the solution~$u$ in
the $\|\cdot\|_H$\nobreakdash-norm, but by the deviation between $\jmath(u)$
and $\jmath(u^h)$ for all $\jmath\in\DD$. Concentrating on the largest deviation,
the {\em support function\/} $s(\DD,u-u^h)$ of~$\DD$ at $u-u^h$ therefore constitutes
an appropriate error measure. The support function $s(\DD,\,\cdot\,):H\to\IR$ of
a closed convex set $\DD\subset{}H^{\star}$ is defined as:
\begin{equation}
\label{eq:errmeas}
s(\DD,v)
=
\sup\big\{\langle\jmath,v\rangle:\jmath\in\DD\big\}
\end{equation}
The supremum in~\EQ{errmeas} is in fact attained by a functional in~$\DD$
by virtue of~$\DD\subset{}H^{\star}$ being closed. A maximizing functional
$\jsup\in{}H^{\star}$ such that $\jmath(v)\leq\jsup(v)$ for all $\jmath\in\DD$ is referred to as a {\em supporting functional\/} of~$\DD$ at~$v$.

\paragraph{Remark}
An equivalent characterization of the worst-case multi-objective error $s(\DD,u-u^h)$
is provided by
\begin{equation*}
\sup\big\{\langle{}r(u^h),z\rangle:z\in{}H\text{ such that }a(v,z)=\jmath(v)\text{ for all }v\in{}H,\:\jmath\in\DD\big\}
\end{equation*}
with $r:H\to{}H^{\star}$, $\langle{}r(u),\,\cdot\,\rangle:=b(\cdot)-a(u,\cdot)$
the residual
functional associated with~\EQ{P} at~$u$; cf. \EQ{residual} and \EQ{dual0}--\EQ{erf}.
Young--Fenchel duality implies that the supremum in~\EQ{supr0} coincides with
the support function of~$\DD$ at~$x\in{}H$, that is
\begin{equation}
\label{eq:supr0}
s(\DD,x) = \sup\big\{\langle{}r(u^h),z\rangle:z\in{}H\text{ such that }a(v,z)=\jmath(v)\text{ for all }v\in{}H,\:\jmath\in\DD\big\}
\end{equation}
provided
\begin{equation}
\label{eq:errorrelation}
x\in{}H:\qquad{}a(x,v)=\langle{}r(u^h),v\rangle\qquad\forall{}v\in{}H
\end{equation}
Noting that $\langle{}r(u^h),v\rangle=b(v)-a(u^h,v)=a(u-u^h,v)$, we infer
that $x=u-u^h$ and, indeed, the supremum in~\EQ{supr0} coincides with $s(\DD,u-u^h)$.
\closeremark

Based on a supporting functional $\jsup\in{}H^{\star}$ of the objective set~$\DD$
at the error~$u-u^h$, an upper bound to the worst-case multi-objective error can be constructed
which depends on the approximation~$u^h$ only via the corresponding residual, $r(u^h)$.
The worst-case multi-objective error bound assumes the conventional DWR form; cf.~\EQ{dual0}
and~\EQ{supr0}.
Denoting by $z\in{}H$ the dual solution associated with the
considered supporting functional of~$\DD$ at~$u-u^h$,
\begin{equation}
\label{eq:duals}
z\in{}H:\qquad{}a(w,z)=\jsup(w)\qquad\forall{}w\in{}H.
\end{equation}
the following sequence of inequalities holds:
\begin{equation}
\label{eq:multiobjectSeq}
\sup_{\jmath\in\DD}\big|\jmath(u)-\jmath(u^h)\big|
\leq
\jsup(u-u^h)
=a(u-u^h,z)
=b(z)-a(u^h,z)=\langle{}r(u^h),z\rangle
\end{equation}
The upper bound in~\EQ{multiobjectSeq} is sharp if $\jsup\in{}\DD$. If the approximation~$u^h$
corresponds to the Galerkin approximation in~\EQ{Ph}, then $\langle{}r(u^h),v^h\rangle$ vanishes for all $v^h\in{}H^h$ and, accordingly, any interpolant of~$z$ in~$H^h$
can be subtracted in the ultimate expression in~\EQ{multiobjectSeq}.

Application of the above theory in actual computations requires the construction of (a suitable approximation of) the support function of the objective set~$\DD$ at the error
$u-u^h$. In addition, if the worst-case multi-objective error estimate serves to direct an adaptivity procedure based on a DWR formulation, then the supporting functional of~$\DD$ at $u-u^h$ is also required.
Below we consider support functions and supporting functionals for three
specific classes of~$\DD$. Exemplifications in the context of the model problem of Section~\SEC{ModProb} are provided in Section~\SEC{numer-exper}.

\paragraph{Example 1} In our first example we consider the case of~$\DD$ represented by a
convex combination of a finite set of functionals, i.e. given $\{\jmath_1,\ldots,\jmath_n\}\subset{}H^{\star}$,
\begin{equation}
\DD=\bigg\{\jmath\in{}H^{\star}:\jmath=\sum_{i=1}^n\alpha_i\jmath_i,\,\sum_{i=1}^n\alpha_i=1,\,\alpha_i\in\IR_{\geq0}\bigg\}
\end{equation}
For any error $u-u^h\in{}H$, it then follows that there exists an index $\check{\imath}\in[1,\ldots,n]$ such that $\jmath_{\check{\imath}}$ is a supporting functional of~$\DD$
at~$u-u^h$:
\begin{multline}
s(\DD,u-u^h)
=
\sup\big\{\langle\jmath,u-u^h\rangle:\jmath\in\DD\big\}
=
\sup\bigg\{\sum_{i=1}^n\alpha_i\langle\jmath_i,u-u^h\rangle: \sum_{i=1}^n\alpha_i=1,\,\alpha_i\in\IR_{\geq0}\bigg\}
\\
=
\max\big\{\langle\jmath_i,u-u^h\rangle:i\in[1,\ldots,n]\big\}
=
\langle\jmath_{\check{\imath}},u-u^h\rangle
\end{multline}
Note that the index~$\check{\imath}$ generally depends on $u-u^h$. The worst-case multi-objective error
$s(\DD,u-u^h)$ can thus be expressed as $s(\DD,u-u^h) = \langle\jmath_{\check{\imath}},x\rangle$ with~$x$ according to~\EQ{errorrelation}.

\paragraph{Example 2} Our second example concerns the case that~$\DD$ coincides with the unit ball in~$H^{\star}$:
\begin{equation}
\label{eq:O=Hneq}
\DD=\{\jmath\in{}H^{\star}: \|\jmath\|_{H^\star} \le 1\}\,.
\end{equation}
For any error $u-u^h$, a supporting functional of~$\DD$ according to~\EQ{O=Hneq}
at~$u-u^h$ is then provided by~$\smash[tb]{(u-u^h,\cdot)_H/\|u-u^h\|_H}$.
Indeed, denoting by $\phi:H^{\star}\to{}H$ the canonical isometry between~$H^{\star}$
and~$H$, it holds for any $u-u^h\in{}H$ that
\begin{multline}
s(\DD,u-u^h)
=
\sup\big\{\langle\jmath,u-u^h\rangle:\jmath\in{}H^{\star},\|\jmath\|_{H^{\star}}\leq1\big\}
=
\sup\big\{
(\phi(\jmath),u-u^h)_{H}:\jmath\in{}H^{\star},\|\jmath\|_{H^{\star}}\leq1\big\}
\\
=
\sup\big\{
(v,u-u^h)_{H}:v\in{}H,\|v\|_{H}\leq1\big\}
=(u-u^h,u-u^h)_H/\|u-u^h\|_H
\end{multline}
Hence, $s(\DD,u-u^h)=(u-u^h,x)/\|u-u^h\|_H$ with $x$ according to~\EQ{errorrelation}.

\paragraph{Example 3} As a preliminary to describe the third example, we consider a Hilbert
space~$V\supseteq{}H$ with inner product $(\cdot,\cdot)_V$ and norm $\|\cdot\|_V$,
such that the embedding $\imath:H\to{}V$ is continuous. Let~$B:V\to{}V$ denote
a continuous linear operator. The third example pertains to the objective set
\begin{equation}
\DD=\Big\{\jmath\in{}H^{\star}:\jmath(\cdot)=\big(\imath(\cdot),Bv\big)_V,v\in{}V,\|v\|_V\leq{}1\Big\}
\end{equation}
That is, $\DD$ corresponds to the $V$-inner product with the image of the unit ball
in $V$ under the linear operator~$B$. For any $u-u^h\in{}H$,
the following sequence of identities holds:
\begin{equation}
\label{eq:SeqId12}
\begin{aligned}
s(\DD,u-u^h)&=\sup\big\{\langle\jmath,u-u^h\rangle:\jmath\in\DD\big\}
\\
&=\sup\big\{\big(\imath(u-u^h),Bv\big)_V:v\in{}V,\|v\|_V\leq{}1\big\}
\\
&=\sup\big\{\big(B^{\star}\imath(u-u^h),v\big)_V:v\in{}V,\|v\|_V\leq{}1\big\}
\\
&=\big(B^{\star}\imath(u-u^h),B^{\star}\imath(u-u^h)\big)_V/\|B^{\star}\imath(u-u^h)\|_V
\end{aligned}
\end{equation}
with $B^{\star}:V\to{}V$ the adjoint operator of~$B$ in~$V$. From~\EQ{SeqId12} we
infer that
\begin{equation}
\label{eq:supfunEx3}
s(\DD,u-u^h)=
\frac{\big(B^{\star}\imath(u-u^h),B^{\star}\imath(x)\big)_V}{\|B^{\star}\imath(u-u^h)\|_V}
\end{equation}
with $x$ according to~\EQ{errorrelation}. A particular instance of this third example
in the context of the model problem of Section~\ref{sec:ProbStat} is addressed in Section~\ref{sec:numer-exper}.

\subsection{Worst-case multi-objective error estimation with data incompatibility}
\label{sec:multi-object-error-incompatible}
To provide a generic setting for worst-case multi-objective error estimation
with data incompatibility, we consider a Hilbert space~$H$ and a proper
subspace~$H_0\subset{}H$. Let $a\in[H\times{}H]^{\star}$ and $b\in{}H^{\star}$ denote
continuous bilinear and linear forms, respectively. In addition, we assume
that~$a$ is coercive on $H_0$, i.e. $\smash[tb]{a(u,u)\geq\underline{c}_a\|u\|_H^2}$ for some
constant $\smash[tb]{\underline{c}_a>0}$ for all $u\in{}H_0$.
Given an element $\ell\in{}H$ and an
objective functional $\jmath\in{}H^{\star}$, we consider the measurement problem
of determining $\jmath(u)$ with $u=\ell+u_0$ and $u_0\in{}H_0$ according to:
\begin{equation}
u_0\in{}H_0:
\qquad
a(u_0,v)=b(v)-a(\ell,v)
\qquad\forall{}v\in{}H_0
\end{equation}
Suppose that $\ell$ is not explicitly available, but instead we only have access to
an approximation $\ell^h\in{}H$ with data error $\ell-\ell^h$.
Given a subspace $\smash[tb]{H_0^h\subset{}H_0}$, we consider the
approximation $\smash[tb]{u^h=\ell^h+u_0^h}$,
\begin{equation}
u_0^h\in{}H^h_0:
\qquad
a(u_0^h,v^h)=b(v^h)-a(\ell^h,v^h)
\qquad
\forall{}v^h\in{}H^h_0
\end{equation}
and, in particular, the corresponding approximate measurement, $\jmath(u^h)$. Without loss of generality, we assume that the objective functional can be written in the
general form $\jmath(\cdot)=a(\cdot,\vartheta)$ for some $\vartheta\in{}H$; cf.~\EQ{jform}.
A derivation analogous to~\EQ{erf} then conveys that
\begin{equation}
\label{eq:erfgen}
\jmath(u)-\jmath(u^h)
=
a(u-u^h,\vartheta)=a(\ell-\ell^h,\vartheta-\pi_0\vartheta)+\langle{}r(u^h),\pi_0\vartheta-v^h\rangle
\end{equation}
for all $v^h\in{}H^h_0$ and with $\pi_0\vartheta\in{}H_0$ the projection of $\vartheta$
onto~$H_0$ according to the dual problem:
\begin{equation}
\label{eq:dualgen}
\pi_0\vartheta\in{}H_0:
\qquad
a(v,\pi_0\vartheta)=a(v,\vartheta)
\qquad
\forall{}v\in{}H_0
\end{equation}
Equation~\EQ{erfgen} imparts that the data error contributes to the measurement
error only via $\vartheta-\pi_0\vartheta$. We can therefore restrict our
analysis of data-incompatibility errors to objective sets $\smash[tb]{\DD\subseteq\{a(\,\cdot\,,\vartheta):\vartheta\in{}H_0^{\bot}\}}$ with $\smash[tb]{H_0^{\bot}:=(1-\pi_0)H}$.

We assume that the objective set $\DD$ is generated by the image of the
map $\lambda\mapsto{}a(\,\cdot\,,L(\lambda))$, where $L:T\to{}H_0^{\bot}$ corresponds to a lineair
operator from a Hilbert space $T$ into $H_0^{\bot}$. The inner
product and norm associated with~$T$ are denoted by $(\cdot,\cdot)_T$
and~$\|\cdot\|_T$, respectively. In particular, we consider the objective set:
\begin{equation}
\label{eq:OsetT}
\DD=\big\{\jmath\in{}H^{\star}:\jmath(\cdot)=a(\,\cdot\,,L({\lambda})),\lambda\in{}T,\|\lambda\|_T\leq{}1\big\}
\end{equation}
By virtue of the continuity of~$L:T\to{}H_0^{\bot}$ and of $a\in[H\times{}H]^{\star}$,
it holds that $\hat{\jmath}(\cdot):=a(u-u^h,L(\cdot))\in{}T^{\star}$.
Denoting by $\phi:T^{\star}\to{}T$ the canonical isometry
between $T^{\star}$ and~$T$, the worst-case multi-objective data-error can
be expressed as:
\begin{equation}
\label{eq:dataerrmo}
\sup\big\{\hat{\jmath}(\lambda):\lambda\in{}T,\|\lambda\|_T\leq{}1\big\}
=
\sup\big\{(\phi(\hat{\jmath}),\lambda)_T:\lambda\in{}T,\|\lambda\|_T\leq{}1\big\}
\end{equation}
The supremum in~\EQ{dataerrmo} is attained at $\lambda=\phi(\hat{\jmath})/\|\phi(\hat{\jmath})\|_T$,
yielding the value $\|\phi(\hat{\jmath})\|_T$.

Instead of the error representation corresponding to
the ultimate expression in~\EQ{erfgen}, we have used the
error representation according to the second expression in~\EQ{erfgen} to formulate
the worst-case multi-objective error representation in~\EQ{dataerrmo}.
The reason for this is that in practice, the lift $\ell$ which appears in the ultimate
expression in~\EQ{erfgen} will not be explicitly available, and needs to be replaced
by a suitable approximation. It is then more convenient and not more restrictive
to apply the second expression in~\EQ{erfgen} with the exact solution~$u$ replaced
by a suitable approximation. Based on various orthogonality properties, many distinct
and equivalent forms of the error representation can be derived;
see Sections~\SEC{numer-exper-multi-object-error-incompatible}
and~\SEC{numer-exper-multi-object-adapt-incompatible}.

\section{Numerical Experiments}
\label{sec:numer-exper}
We consider four different numerical experiments to elucidate the main properties of the
worst-case multi-objective error estimator and its application in adaptive-refinement procedures. The first
test case addresses the properties of the error estimate without data incompatibility
with various polynomial orders of the finite-element approximation, for a problem with regular primal and dual solutions. The second
test case illustrates the application of the error estimate in the absence of data incompatibility
to guide an adaptive-refinement process, for a problem in which the primal and dual solutions exhibit restricted
regularity on account of a corner singularity. The third and final test cases are concerned with error estimation
and adaptivity for a problem with nearly singular incompatible boundary data and heterogeneous coefficients.

\subsection{Worst-case multi-objective error estimation without data incompatibility}
\label{sec:numer-exper-multi-object-error-compatible}
The first test case concerns the Dirichlet--Poisson problem~\EQ{BVP} on a bi-unit square
$\Omega=(-1,1)^2$. The data~$f$ and $g$ are selected such that
\newcommand{\uosc}{u_{\mathrm{osc}}}
\begin{equation}
\label{eq:usample1}
\uosc(x,y)=\sin(6\pi{}x)\sin(6\pi{}y)
\end{equation}
solves~\EQ{BVP}. Note that the solution $\uosc$ is smooth (of class $C^{\infty}$) and displays moderate
oscillatory behavior, and that the trace of $\uosc$
on~$\partial\Omega$ vanishes and, hence, the corresponding Dirichlet data~$g$ is homogeneous.
The considered Dirichlet--Poisson problem can be condensed into the weak formulation:
\begin{equation}
\label{eq:TC1weak}
u\in{}H^1_0(\Omega):\qquad{}a(u,v)=b(v)\qquad\forall{}v\in{}H^1_0(\Omega)
\end{equation}
with the bilinear form $a(\cdot,\cdot)$ corresponding to~\EQ{a(u,v)b(v)} and the linear form $b(\cdot)=a(\uosc,\cdot)$,
with $\uosc$ according to~\EQ{usample1}.
By virtue of the homogeneity of the Dirichlet data, the finite-element approximations are devoid of data-incompatibility errors. Moreover, the smoothness of~$\uosc$ ensures that finite-element approximations on uniform meshes display optimal convergence behavior, in particular, $\|u^h-\uosc\|_{H^k(\Omega)}\leq{}C{}h^{p+1-k}$ ($k\in\{0,1\}$) for
some constant $C>0$ independent of the mesh width $h$, with $p$ the polynomial degree of the finite-element space.

We set $H:=H^1_0(\Omega)$, $V:=L^2(\Omega)$ and consider the objective set
\begin{equation}
\label{eq:TC1O}
\DD=\big\{\jmath\in[H^1(\Omega)]^{\star}:\jmath(\cdot)=(\,\cdot\,,\chi_{\omega}v)_{L^2(\Omega)},\|v\|_{L^2(\Omega)}\leq{}1\big\}
\end{equation}
with $\chi_{\omega}$ the indicator function of the set $\omega=(-3/4,-1/4)^2$.
The support function of $\DD$ at~$u-u^h$ characterizes
the $L^2$-norm of the error in $u^h$ in the subset~$\omega$. Indeed, it holds that
$\chi_{\omega}^{\star}=\chi_{\omega}$ and, hence,
\begin{equation}
\begin{aligned}
\sup\big\{(u-u^h,\chi_{\omega}v)_{L^2(\Omega)}:\|v\|_{L^2(\Omega)}\leq{}1\big\}
&=
\sup\big\{\big(\chi_{\omega}(u-u^h),v\big)_{L^2(\Omega)}:\|v\|_{L^2(\Omega)}\leq{}1\big\}
\\
&=
\big\|\chi_{\omega}(u-u^h)\big\|_{L^2(\Omega)}
\end{aligned}
\end{equation}
The considered Dirichlet--Poisson problem with objective set~$\DD$ according
to~\EQ{TC1O} conforms to the third example treated in Section~\ref{sec:multi-object-error}. By analogy to~\EQ{supfunEx3}, we infer that
\begin{equation*}
\jsup(\cdot)=
\frac{(\chi_{\omega}(u-u^h),\,\cdot\,)_{L^2(\Omega)}}{\|u-u^h\|_{L^2(\Omega)}}
\end{equation*}
is a supporting functional of~$\DD$
at $u-u^h$. Therefore, if $u^h$ is determined from the Galerkin approximation~\EQ{GalerkinPrimal}, then $s(\DD,u-u^h)=\langle{}r(u^h),z-v^h\rangle$ for any $v^h\in{}V^h_0$, with $z\in{}H^1_0(\Omega)$ according to the dual problem
\begin{equation}
\label{eq:dualTC1}
z\in{}H^1_0(\Omega):\qquad a(v,z)=\frac{(\chi_{\omega}(u-u^h),v)_{L^2(\Omega)}}{\|u-u^h\|_{L^2(\Omega)}}\qquad\forall{}v\in{}H^1_0(\Omega)
\end{equation}

In practice, the exact solution $u$ is unavailable, and an estimate of the error $u-u^h$
is required to construct an approximate supporting functional. In the numerical
computations, we apply finite-element approximation spaces based on rectangular elements,
and we opt to approximate the error in~\EQ{dualTC1} by $u^{h/2}-u^h$, where~$u^{h/2}$ is an approximation obtained from a Galerkin approximation in a finite-element
space $\smash[tb]{V^{h/2}_0}$ in which each element is uniformly divided into 4 elements; see Figure~\FIG{WCMO_Alg}. The refined finite-element space  $\smash[tb]{V^{h/2}_0}$ moreover serves to construct an approximation to the dual solution in~\EQ{dualTC1}. In summary, considering an
approximation $\smash[tb]{u^h\in{}V^h_0}$ and a refined approximation $\smash[tb]{u^{h/2}\in{}V^{h/2}_0}$,
the worst-case multi-objective error estimate pertaining to $u^h$
is $\smash[tb]{\langle{}r(u^h),z^{h/2}\rangle}$ with $\smash[tb]{z^{h/2}}$ according to
\begin{equation}
\label{eq:dualhTC1}
z^{h/2}\in{}V^{h/2}_0:\qquad{}a(v^{h/2},z^{h/2})=
\frac{(\chi_{\omega}(u^{h/2}-u^h),v^{h/2})_{L^2(\Omega)}}{\|\chi_{\omega}(u^{h/2}-u^h)\|_{L^2(\Omega)}}
\qquad\forall{}v^{h/2}\in{}V^{h/2}_0
\end{equation}

\begin{figure}[t]
\begin{center}
\includegraphics[width=0.9\textwidth]{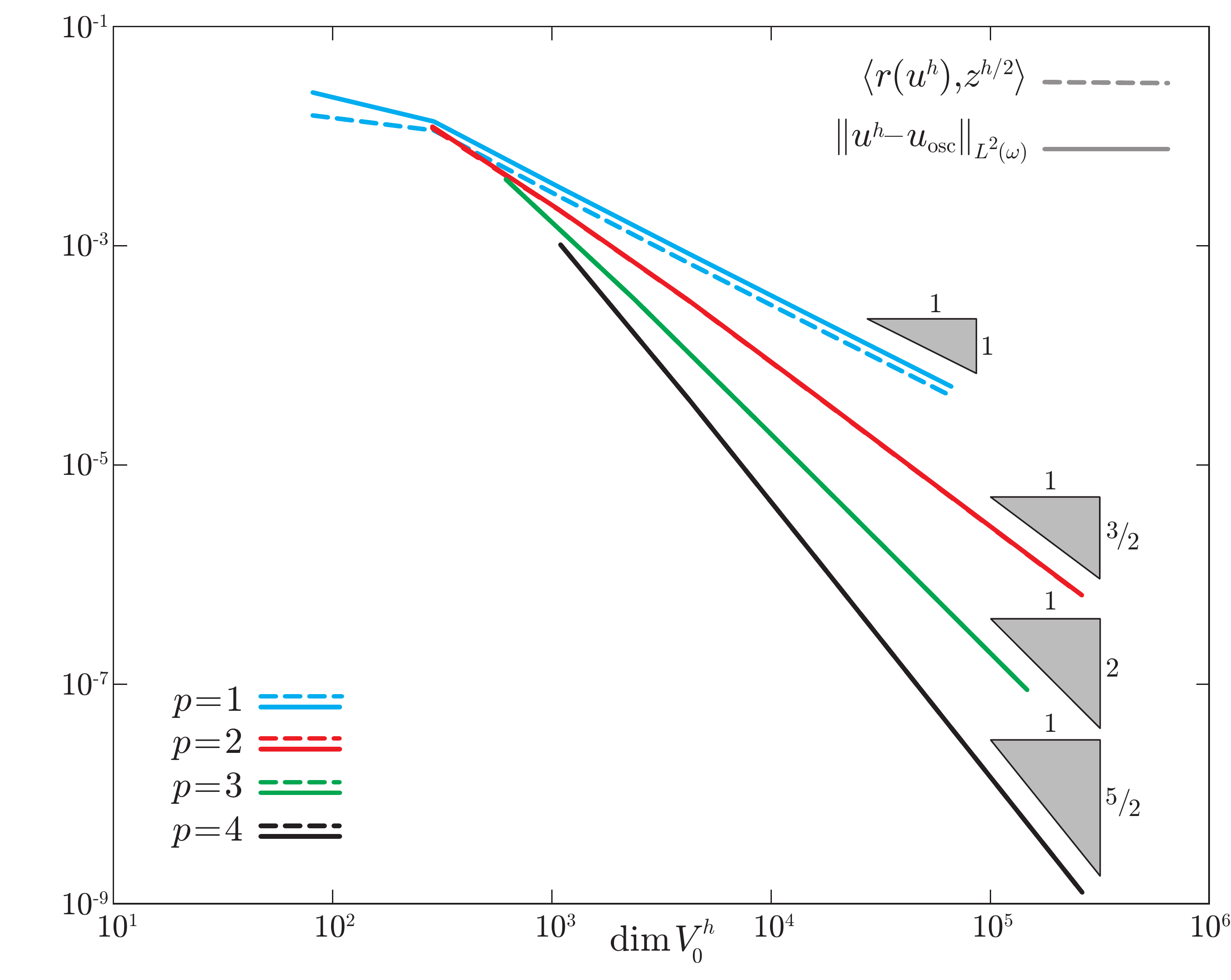}
\caption{%
Worst-case multi-objective error estimate
$\smash[tb]{\langle{}r(u^h),z^{h/2}\rangle}$
for the Dirichlet--Poisson problem with solution~\EQ{usample1} and objective set~\EQ{TC1O} versus $\smash[tb]{\dim{}V^h_0}$
for uniform finite-element spaces~$\smash[tb]{V^h_0}$ with polynomial orders $p\in\{1,2,3,4\}$. The solid lines indicate the actual errors $\|u^h-\uosc\|_{L^2(\omega)}$.
\label{fig:TC1_conv}}
\end{center}
\end{figure}

Figure~\FIG{TC1_conv} plots the worst-case multi-objective error estimate, $\langle{}r(u^h),z^{h/2}\rangle$,
versus the dimension of the finite-element approximation space, $\smash[tb]{\dim{}V^h_0}$,
for finite-element approximations with polynomial degrees $p\in\{1,2,3,4\}$ on uniform meshes with mesh width
$h=2^{-2},2^{-3},\ldots$. To verify the accuracy of the worst-case multi-objective error estimate, the figure
also plots the exact error $\|u^h-\uosc\|_{L^2(\omega)}$. It is noteworthy that both the
worst-case multi-objective error estimate and the exact error exhibit optimal convergence rates proportional
to~$(\smash[tb]{\dim{}V^h_0})^{-(p+1)/2}$. The deviation between the worst-case multi-objective error estimate and the exact error is generally small. For $p=1$, the error estimate underestimates the error by approximately 25\% on the
coarsest mesh and by approximately 10\% on all finer meshes. For $p\in\{2,3,4\}$, the deviation between the error estimate
and the exact error is not discernible.

\subsection{Worst-case multi-objective error estimation and adaptivity without data incompatibility}
\label{sec:numer-exper-multi-object-adapt-compatible}
\newcommand{\using}{u_{\mathrm{sing}}}
In the second test case we examine the application of the worst-case multi-objective error estimate
to direct an adaptive-refinement procedure, following the standard SEMR
(\texttt{Solve} $\rightarrow$ \texttt{Estimate} $\rightarrow$ \texttt{Mark}
$\rightarrow$ \texttt{Refine}) process; see, for instance, \cite{Nochetto:2012hl,Dorfler:1996uq}.
We consider the Dirichlet--Poisson problem~\EQ{BVP} on an L\nobreakdash-shaped domain $\Omega=(-1,1)^2\setminus([0,1)\times(-1,0])$.
The data~$f$ and $g$ are selected such that $u=\using+\uosc$ solves~\EQ{BVP} with
\begin{equation}
\label{eq:using}
\using{}(x,y)=\cos\bigg(\frac{\pi{}x}{2}\bigg)\cos\bigg(\frac{\pi{}y}{2}\bigg)\big(x^2+y^2\big)^{1/3}
\sin\bigg(\frac{2}{3}\arctan\Big(\frac{y}{x}\Big)\bigg)
\end{equation}
and $\uosc$ according to~\EQ{usample1}. Note that $\gamma{}u=0$ and, hence, the corresponding
Dirichlet data $g$ is homogeneous.
The function~$\using$ in~\EQ{using} exhibits a corner singularity at $(0,0)$. In particular, it holds
that $\using\in{}H^s(\Omega)$ for $s<5/3$; see, for instance, \cite{Nochetto:2012hl}.
It is well known that finite-element approximations on uniform meshes display suboptimal convergence rates,
owing to the singularity in the solution. On uniform meshes, the error $\|u^h-\using\|_{H^1(\Omega)}$
in the best approximation in $\smash[tb]{V^h_0}$ is only proportional to $\smash[tb]{h^{2/3}}$, independent
of the order of approximation~\cite{Nochetto:2012hl}.
Moreover, an Aubin-Nitsche duality argument conveys that $\smash[tb]{\|u^h-\using\|_{L^2(\Omega)}}$ decays
only as~$\smash[tb]{h^{4/3}}$. The purpose of an adaptive finite-element method is to restore
optimal convergence rates with respect to the dimension of the approximation space.

\begin{figure}
\begin{center}
\centering\epsfig{figure=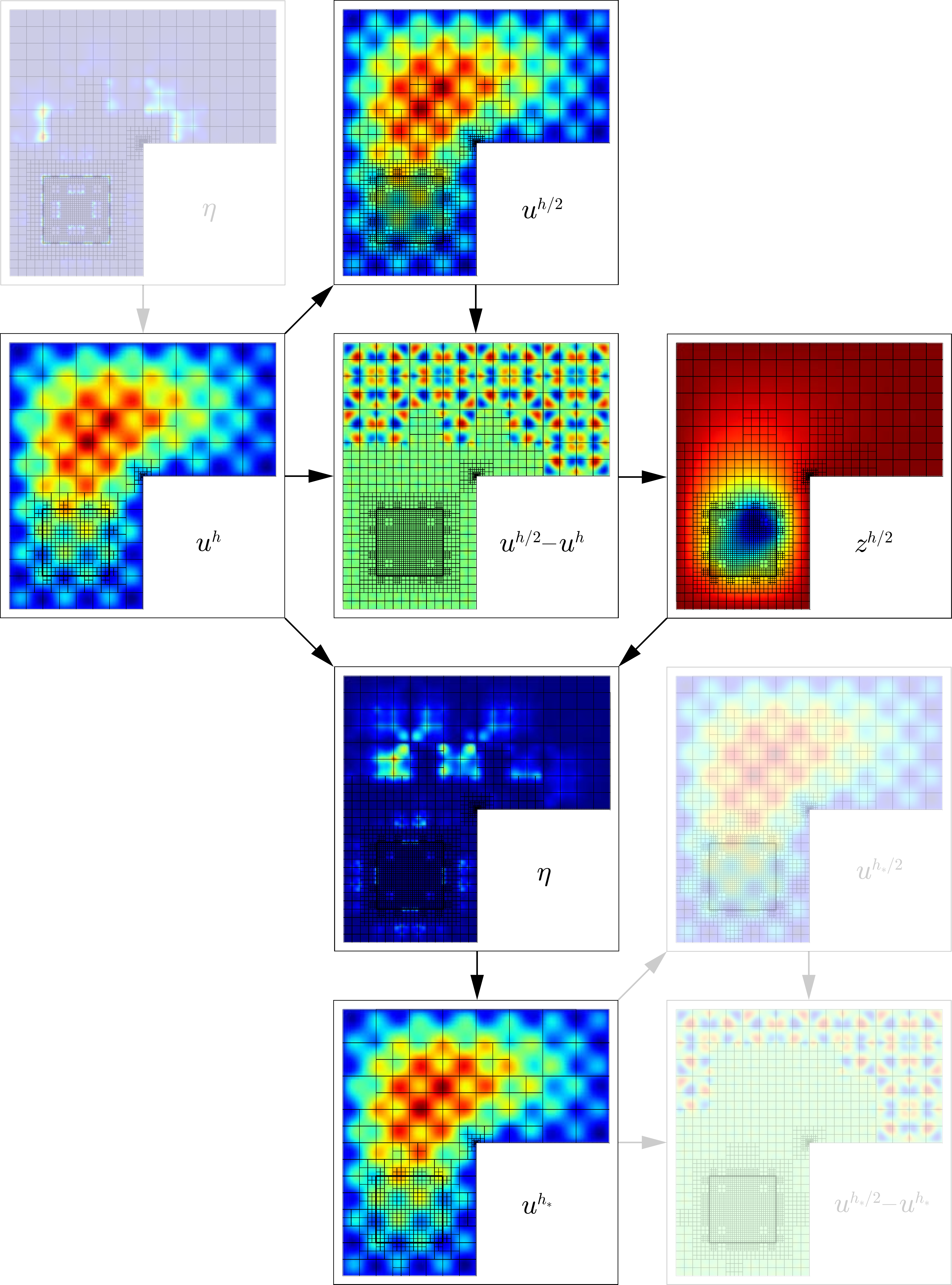,width=0.94\textwidth}
\caption{%
Illustration of the worst-case multi-objective error-estimation process and the corresponding adaptive-refinement procedure.
The square in each panel indicates the set $\omega=(-3/4,-1/4)^2$ in which the $L^2$-norm of the error
is to be minimized, in accordance with the objective set $\mathcal{O}$ in~\EQ{TC1O}.
\label{fig:WCMO_Alg}}
\end{center}
\end{figure}

The objective set for the worst-case multi-objective error-estimation procedure is again selected
in accordance with~\EQ{TC1O} and, hence, the support function of~$\mathcal{O}$ corresponds to the
$\|\cdot\|_{L^2(\omega)}$\nobreakdash-norm of the
error in the finite element approximation in the set $\omega=(-3/4,-1/4)^2$; see Figure~\FIG{WCMO_Alg} for an illustration.
The adaptive-refinement process based on the aforementioned error estimate thus endeavors to restore the optimal
rate of convergence of $\|u^h-u\|_{L^2(\omega)}$ with respect to $\dim{}\smash[tb]{V^h_0}$. The  applied SEMR process
is illustrated in Figure~\FIG{WCMO_Alg}. Given a mesh $\mathcal{T}^h$, we determine the finite-element approximation
in the corresponding finite-element space, $u^h\in\smash[tb]{V^h_0}$; see panel~(2,1) of Figure~\FIG{WCMO_Alg}.
We then construct the refined mesh, $\smash[tb]{\mathcal{T}^{h/2}}$, by uniformly subdividing each element of~$\smash[tb]{\mathcal{T}^h}$ into~$4$ elements, and determine the finite-element approximation $u^{h/2}\in\smash[tb]{V^{h/2}_0}$; see panel~(1,2). Based on $u^{h/2}$ and $u^h$, the error in the finite-element solution in~$\smash[tb]{V^h_0}$
is estimated as $u^{h/2}-u^h$ (panel~(2,2)). The error estimate $u^{h/2}-u^h$ serves to construct the approximate
dual solution $z^{h/2}\in\smash[tb]{V^{h/2}_0}$ corresponding to the supporting functional of the objective set~$\mathcal{O}$ in~\EQ{TC1O}, in accordance with~\EQ{dualhTC1} (panel~(2,3)). Based on $u^h$ and $z^{h/2}$, we determine
the worst-case multi-objective error estimate $\langle{}r(u^h),z^{h/2}\rangle$. The error estimate is subsequently
decomposed into local contributions, associated with the support of basis functions in~$\smash[tb]{V^{h/2}_0}$;
see panel~(3,2). Let us note that instead of a traditional element-wise marking strategy, we apply the function-support
marking strategy introduced in~\cite{Kuru:2014fk} (see also~\cite{Richter:2015zr}). To elaborate on this
strategy, let $\smash[tb]{\{\psi_i\}}_{i\in\II}$
with $\II\subset\mathbb{N}$ denote the applied basis
of~$\smash[tb]{V^{h/2}_0}$ and let $\Pi$ denote a suitable projection
from~$\smash[tb]{V_0^{h/2}}$ to~$\smash[tb]{V_0^{h}}$. In particular, we apply
an $L^2$\nobreakdash-projection for reasons of implementational convenience. The
difference between the approximate dual solution and its projection
onto~$\smash[tb]{V_0^h}$ resides in $\smash[tb]{V^{h/2}_0}$ and can therefore
be expanded as $z^{h/2}-\Pi{}z^{h/2}=\sum\sigma_i\psi_i$. Hence, we can decompose
the DWR error estimate into fine-mesh basis-function contributions as:
\begin{equation}
\label{eq:errorbound1}
\big|\big\langle{}r(u^h),z^{h/2}-\Pi{}z^{h/2}\big\rangle\big|
=
\bigg|\sum_{i\in\II}\sigma_i\langle{}r(u^h),\psi_i\rangle\bigg|
\leq
\sum_{i\in\II}\eta_i
\end{equation}
with $\eta_i=|\sigma_i\langle{}r(u^h),\psi_i\rangle|$. Based on the basis-function
indicators,~$\{\eta_i\}$, we mark basis functions according to
a D\"orfler-type~\cite{Dorfler:1996uq} marking, i.e. we select a minimal
set of indices $\IA\subset\II$ such that
\begin{equation}
\sum_{i\in\IA}\eta_i\geq(1-\zeta)\sum_{i\in\II}\eta_i
\end{equation}
for some $\zeta\in(0,1)$. In particular, we set $\zeta=1/2$. Finally, we select a
minimal set of coarse-grid elements, $\IM$, such that the union of the supports of
all marked basis functions, $\cup_{i\in\IA}\mathrm{supp}(\psi_i)$, is contained in the union of the selected elements, $\cup_{K\in\IM}K$. These elements are subsequently
partitioned to form a new mesh and to construct a corresponding new refined finite-element space; see panel (4,2).

\begin{figure}
\begin{center}
\centering\epsfig{figure=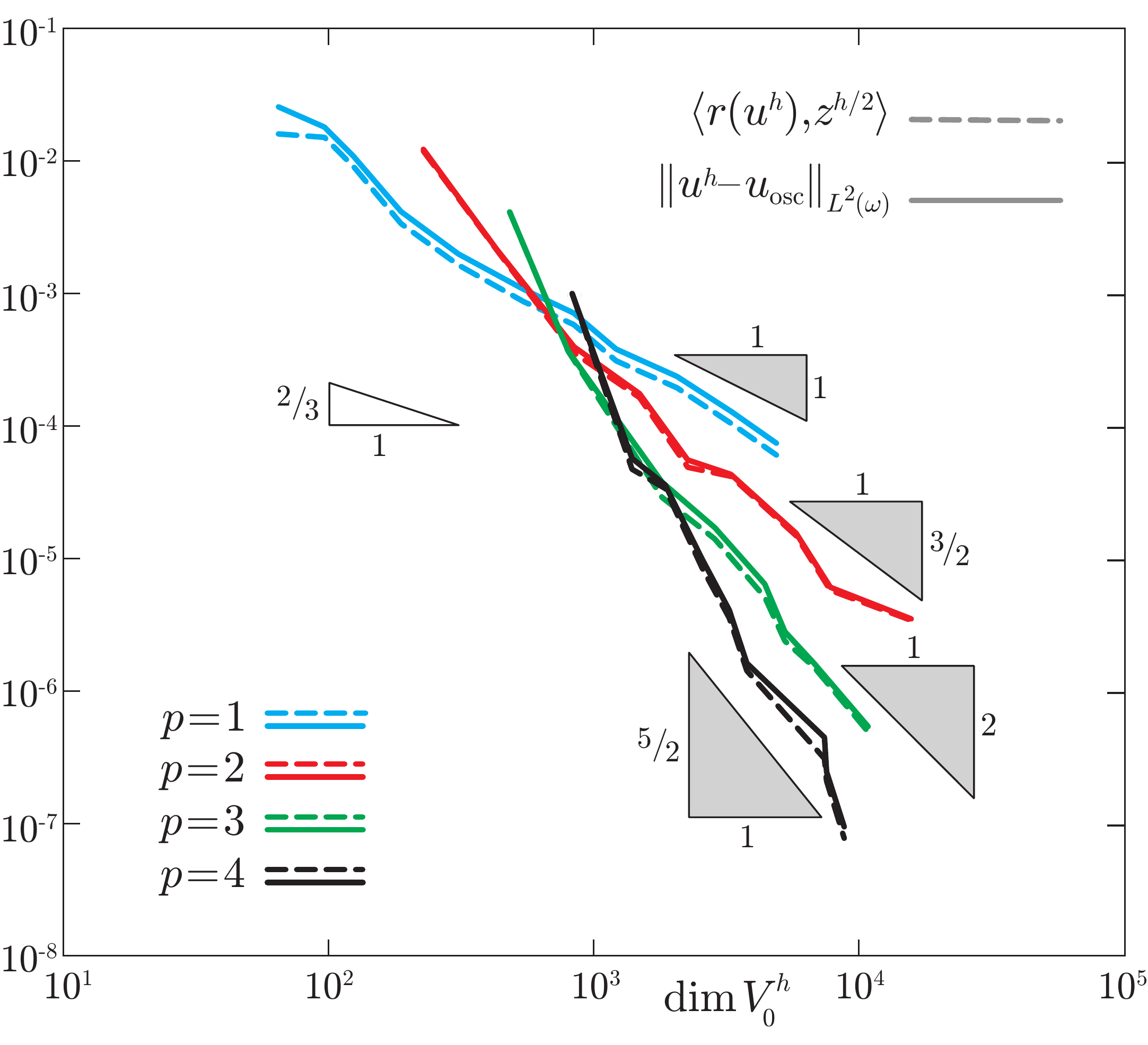,width=0.9\textwidth}
\caption{%
Worst-case multi-objective error estimate $\smash[tb]{\langle{}r(u^h),z^{h/2}\rangle}$ 
and exact error $\|u^h-u\|_{L^2(\omega)}$ versus $\smash[tb]{\dim{}V^h_0}$,
for sequences of adaptive finite-element approximations with polynomial orders $p\in\{1,2,3,4\}$
of the Dirichlet--Poisson problem on an L\nobreakdash-shaped domain with solution $u=\uosc+\using$. The white triangle
indicates the ($p$\nobreakdash-independent) convergence rate corresponding to uniform mesh refinement.
\label{fig:TC2_conv}}
\end{center}
\end{figure}

Figure~\FIG{TC2_conv} plots the worst-case multi-objective error estimate $\langle{}r(u^h),z^{h/2}\rangle$ and
the error $\|u^h-u\|_{L^2(\omega)}$ corresponding to the finite-element approximation~$u^h$ generated by the SEMR algorithm, versus the dimension of the adaptively refined finite-element space, $\smash[tb]{\dim{}V^h_0}$,
for finite-element approximations with polynomial orders $p\in\{1,2,3,4\}$. The results convey that for all
considered polynomial orders, the SEMR algorithm restores the optimal rate of convergence, i.e.,
the sequence of finite-element approximations generated by the SEMR algorithm exhibits a convergence rate proportional
to $\smash[tb]{(\dim{}V^h_0)^{-(p+1)/2}}$. It is noteworthy that the deviation between the actual error and the worst-case
multi-objective error estimate is generally small.

\subsection{Worst-case multi-objective error estimation with data incompatibility}
\label{sec:numer-exper-multi-object-error-incompatible}
\begin{figure}
\begin{center}
\centering\epsfig{figure=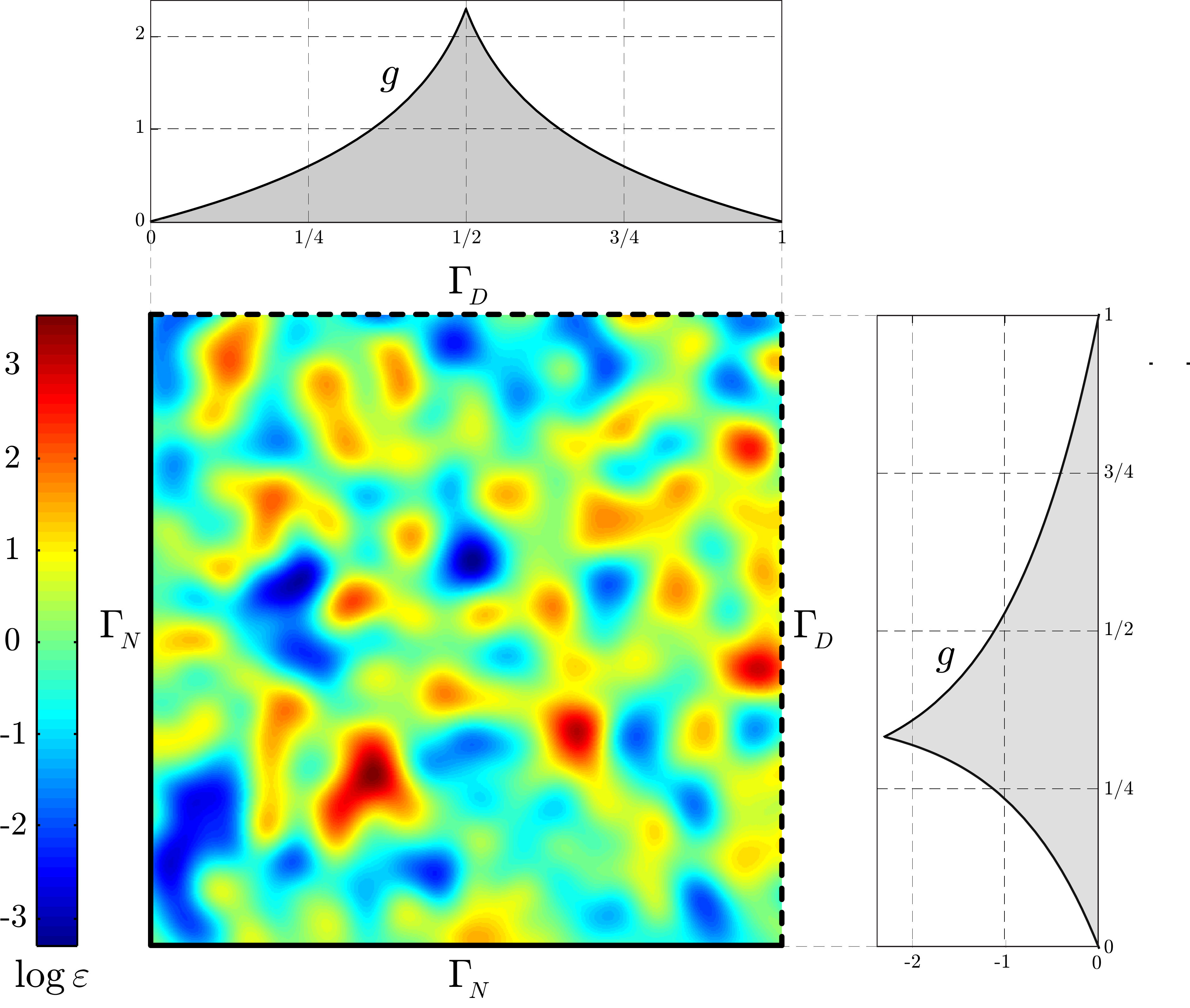,width=0.90\textwidth}
\caption{%
Illustration of the setup of the third test case. Colors encode the logarithm of the
heterogeneous permeability, $\log\varepsilon$.
\label{fig:TC3_Setup}}
\end{center}
\end{figure}

To illustrate the worst-case multi-objective error-estimation procedure with data incompatibility,
we consider the unit square $\Omega=(0,1)^2$ and the Dirichlet--Neumann--Laplace problem:
\begin{subequations}
\label{eq:BVPTC3}
\begin{alignat}{2}
-\grad\cdot(\varepsilon\grad{}u)&=0&\quad&\text{in }\Omega
\label{eq:DiffEqTC3}
\\
u&=g&\quad&\text{on }\Gamma_D
\label{eq:BCsDTC3}
\\
\varepsilon\partial_nu&=0&\quad&\text{on }\Gamma_N
\label{eq:BCsNTC3}
\end{alignat}
\end{subequations}
with $\Gamma_D=((0,1]\times\{1\})\cup(\{1\}\times(0,1])$ and
$\Gamma_N=\partial\Omega\setminus\Gamma_D$. The Dirichlet data $g:\Gamma_D\to\IR$
are given by:
\begin{equation}
g=
\begin{cases}
-\log\big((1-2x)+\nu\,2x\big)
&\text{if }(x,y)\in(0,1/2]\times\{1\}
\\
-\log\big((2x-1)+\nu(2-2x)\big)
&\text{if }(x,y)\in(1/2,1]\times\{1\}
\\
\phantom{-}\log\big((1-3x)+\nu\,3x\big)
&\text{if }(x,y)\in\{1\}\times(0,1/3]
\\
\phantom{-}\log\big((3x-1)/2+\nu\,3(1-x)/2\big)
&\text{if }(x,y)\in\{1\}\times(1/3,1]
\end{cases}
\end{equation}
with $\nu=10^{-1}$; see Figure~\FIG{TC3_Setup}. The boundary data has been selected such
that it exhibits nearly singular behavior at $(x,y)\in\{1/2\}\times\{1\}$
and at $(x,y)\in\{1\}\times\{1/3\}$. The location of the latter point is such
that it does not correspond to a vertex of the finite-element mesh for any finite refinement
by uniform subdivision. We consider~\EQ{BVPTC3} with both
a homogeneous coefficient~$\varepsilon=1$ and a heterogeneous coefficient
as depicted in~\FIG{TC3_Setup}.

Problem~\EQ{BVPTC3} can be conceived of as a Darcy-type problem, with~$u$ acting as pressure
and~$-\varepsilon\grad{}u$ the flux induced by the pressure gradient in a medium with 
isotropic permeability $\varepsilon:\Omega\to\IR_{>0}$. In this setting, the Dirichlet
condition~\EQ{BCsDTC3} imposes a prescribed pressure on~$\Gamma_D$ and the homogeneous Neumann
condition~\EQ{BCsNTC3} corresponds to an impermeability condition on~$\Gamma_N$.
The heterogeneous permeability field is a smoothed version of a typical
permeability field used in geostatistics in hydraulic diffusivity inversion.
The permeability data is provided on a uniform mesh of $628\times628$ points and is linearly
interpolated between the data points to construct a continuous permeability field.

Denoting by $\smash[tb]{H^1_{0,\Gamma_D}(\Omega)}$ the subspace of $H^1(\Omega)$ of functions that
vanish on $\Gamma_D$ in the trace sense, and by $\ell_{g}$ a lift of~$g$
into~$H^1(\Omega)$ such that $\smash[tb]{(\gamma\ell_{g})|_{\Gamma_D}=g}$,
the boundary value problem~\EQ{BVPTC3} can
be condensed into the weak formulation: Determine $u=u_0+\ell_{g}$ with
\begin{equation}
\label{eq:weakForm}
u_0\in{}H^1_{0,\Gamma_D}(\Omega):\qquad{}a_{\varepsilon}(u_0,v)=-a_{\varepsilon}(\ell_{g},v)
\qquad\forall{}v\in{}H^1_{0,\Gamma_D}
\end{equation}
and $a_{\varepsilon}:H^1(\Omega)\times{}H^1(\Omega)\to\IR$ according to:
\begin{equation}
\label{eq:aeps}
a_{\varepsilon}(u,v)=\int_{\Omega}\varepsilon\grad{}u\cdot\grad{}v
\end{equation}

In conjunction with~\EQ{weakForm},
we consider an objective set that is associated with the flux functional on a subset~$\omega\subseteq\Gamma_D$:
\begin{equation}
\label{eq:TC3O}
\DD=\Big\{\jmath\in{}[H^1(\Omega)]^{\star}:\jmath(\cdot)
=
a_{\varepsilon}(\cdot,l_{\lambda}),\lambda\in{}H^{1/2}_{0,\partial\Omega\setminus\omega}(\partial\Omega),\varnorm{\lambda}_{H^{1/2}(\partial\Omega)}\leq{}1\Big\}
\end{equation}
where $H^{1/2}_{0,\partial\Omega\setminus\omega}(\partial\Omega)$ represents the subset of functions in $H^{1/2}(\partial\Omega)$ that vanish on $\partial\Omega\setminus\omega$,
and $l_{\lambda}$ is the Moore--Penrose lift of $\lambda$ with respect to the norm induced
by $a_{\varepsilon}$:
\begin{equation}
\label{eq:MPliftTC3}
l_{\lambda}=\argmin\big\{a_{\varepsilon}(v,v):v\in{}H^1(\Omega),\gamma{}v=\lambda\big\}
\end{equation}
Furthermore, in~\EQ{TC3O} we have introduced the 
norm $\varnorm{\lambda}_{H^{1/2}(\partial\Omega)}^2=\smash[tb]{(\lambda,\lambda)_{H^{1/2}(\partial\Omega)}}$
corresponding to the inner product $\smash[tb]{(\varphi,\lambda)_{H^{1/2}(\partial\Omega)}}
=a_{\varepsilon}(l_{\varphi},l_{\lambda})$. 

The objective set~\EQ{TC3O}
conforms to~\EQ{OsetT} with, in particular, $T=\smash[tb]{H{}^{1/2}_{0,\partial\Omega\setminus\omega}(\partial\Omega)}$ and~$L:\lambda\mapsto{}l_{\lambda}$. 
%
From~\EQ{dataerrmo} we then infer that the worst-case multi-objective error satisfies
\begin{multline}
\label{eq:TC3ids}
\sup\Big\{a_{\varepsilon}(u-u^h,l_{\lambda}):\lambda\in{}H^{1/2}_{0,\partial\Omega\setminus\omega}(\partial\Omega),\varnorm{\lambda}_{H^{1/2}(\partial\Omega)}\leq{}1\Big\}
\\=
\sup\big\{a_{\varepsilon}(l_{\varphi},l_{\lambda}):\lambda\in{}H^{1/2}_{0,\partial\Omega\setminus\omega}(\partial\Omega),\varnorm{\lambda}_{H^{1/2}(\partial\Omega)}\leq{}1\Big\}
=\sqrt{a_{\varepsilon}(l_{\varphi},l_{\varphi})}
\end{multline}
with $\varphi\in{}H^{1/2}_{0,\partial\Omega\setminus\omega}(\partial\Omega)$ according to:
\begin{equation}
\label{eq:varphi1}
\varphi\in{}H^{1/2}_{0,\partial\Omega\setminus\omega}(\partial\Omega):
\qquad
a_{\varepsilon}(l_{\varphi},l_{\lambda})=a_{\varepsilon}(u-u^h,l_{\lambda})
\qquad
\forall{\lambda}\in{}H^{1/2}_{0,\partial\Omega\setminus\omega}(\partial\Omega)
\end{equation}
To elucidate the relation between~\EQ{TC3ids} and~\EQ{dataerrmo}, we note 
that the left member of~\EQ{varphi1} coincides with the inner product 
$\smash[tb]{(\varphi,\lambda)_{H^{1/2}(\partial\Omega)}}$ and that the right member of~\EQ{varphi1}
corresponds to the functional $\hat{\jmath}(\cdot):=a_{\varepsilon}(u-u^h,l_{(\cdot)})\in
\smash[tb]{[H^{1/2}_{0,\partial\Omega\setminus\omega}(\partial\Omega)]^{\star}}$. 
Equation~\EQ{varphi1} therefore constitutes a map $\hat{\jmath}\mapsto\phi(\hat{\jmath})=\varphi$,
corresponding to the canonical isometry:
\begin{equation}
\phi:
[H^{1/2}_{0,\partial\Omega\setminus\omega}(\partial\Omega)]^{\star}\to
H^{1/2}_{0,\partial\Omega\setminus\omega}(\partial\Omega)\\
\end{equation}
Furthermore, it holds that 
$\smash[tb]{\varnorm{\phi(\hat{\jmath})}_{H^{1/2}(\partial\Omega)}^2}=  
\smash[tb]{(\varphi,\varphi)_{H^{1/2}(\partial\Omega)}}=
a_{\varepsilon}(l_{\varphi},l_{\varphi})$.

In computations, it is impractical to extract~$\varphi$ directly from~\EQ{varphi1},
on account of the implicit dependence of the inner product $a_{\varepsilon}(l_{(\cdot)},l_{(\cdot)})$ on the lift
operator~$l_{(\cdot)}$. To reformulate~\EQ{varphi1} into
a tractable equivalent form, we note that the optimality conditions associated
with~\EQ{MPliftTC3} imply that
\begin{equation}
\label{eq:qharm}
a_{\varepsilon}(v,l_{\lambda})=0\qquad \forall v\in{}H^1_0(\Omega)
\end{equation}
By virtue of~\EQ{qharm} and the symmetry of $a_{\varepsilon}$,
Equation~\EQ{varphi1} implies
\begin{equation}
\label{eq:varphi2}
a_{\varepsilon}(l_{\varphi},l_{\lambda}+v)=a_{\varepsilon}(u-u^h+q,l_{\lambda}+v)
\qquad\forall{}v\in{}H^1_0(\Omega)
\end{equation}
provided that $q\in{}H^1_0(\Omega)$ satisfies
\begin{equation}
\label{eq:qeq}
q\in{}H^1_0(\Omega):
\qquad
a_{\varepsilon}(q,v)=-a_{\varepsilon}(u-u^h,v)
\qquad\forall{}v\in{}H^1_0(\Omega)
\end{equation}
Because $(\lambda,v)\mapsto{}l_{\lambda}+v$ provides a bijection between $\smash[tb]{H^{1/2}_{0,\partial\Omega\setminus\omega}(\partial\Omega)\times{}H^1_0(\Omega)}$ and $\smash[tb]{H^1_{0,\partial\Omega\setminus\omega}(\Omega)}$, Equation~\EQ{varphi2} implies that
$l_{\varphi}$ is the unique element of $\smash[tb]{H^1_{0,\partial\Omega\setminus\omega}(\Omega)}$ in compliance with:
\begin{equation}
\label{eq:varphi}
l_{\varphi}\in{}H^1_{0,\partial\Omega\setminus\omega}(\Omega):
\qquad
a_{\varepsilon}(l_{\varphi},v)
=
a_{\varepsilon}(u-u^h+q,v)
\qquad\forall{}v\in{}H^1_{0,\partial\Omega\setminus\omega}(\Omega)
\end{equation}
Problem~\EQ{varphi} is in canonical form. Hence, $l_{\varphi}$ can be
conveniently extracted from~\EQ{varphi} and $\varphi$ can in turn be obtained
directly from the identity $\varphi=\gamma{}(l_{\varphi})$.

The exact solution~$u$ is generally not available. Also, Eqs.~\EQ{qeq}\nobreakdash-\EQ{varphi} cannot generally be solved. In the numerical computations, we therefore have to replace~$u$ and~$q$ by
computable approximations. We opt to replace $u$ by an approximation $\smash[tb]{u^{h/2}}$
obtained in the refined finite-element space $V^{h/2}$. In addition, we
replace $q$ and $l_{\varphi}$ by $q^{h/2}$ and $\smash[tb]{l_{\varphi}^{h/2}}$, respectively,
according to:
\begin{alignat}{3}
q^{h/2}\in{}V^{h/2}_0:&&\qquad{}a_{\varepsilon}(q^{h/2},v^{h/2})&=-a_{\varepsilon}(u^{h/2}-u^h,v^{h/2})&\qquad&\forall{}v^{h/2}\in{}V^{h/2}_0
\label{eq:qh/2}
\\
l_{\varphi}^{h/2}\in{}V^{h/2}_{0,\partial\Omega\setminus\omega}:&&\qquad{}a_{\varepsilon}(l_{\varphi}^{h/2},v^{h/2})&=a_{\varepsilon}(u^{h/2}-u^h+q^{h/2},v^{h/2})&\qquad&\forall{}v^{h/2}\in{}V^{h/2}_{0,\partial\Omega\setminus\omega}
\label{eq:lvarphih/2}
\end{alignat}
with $V^{h/2}_{0,\partial\Omega\setminus\omega}:=V^{h/2}\cap{}H^1_{0,\partial\Omega\setminus\omega}(\Omega)$.

If the worst-case multi-objective error estimate is to be applied in an adaptive-refinement process, then a decomposition of the error into basis-function contributions must be introduced, cf. Equation~\EQ{errorbound1}. Denoting by $\{\psi_i\}$ a basis of~$\smash[tb]{V^{h/2}_{0,\partial\Omega\setminus\omega}}$, there exist coefficients $\smash[tb]{\{\sigma_i\}}$ such that:
\begin{equation}
\frac{l_{\varphi}^{h/2}}{\Varnorm{l_{\varphi}^{h/2}}_{H^{1/2}(\partial\Omega)}}
=\sum\sigma_i\psi_i
\end{equation}
Hence, the worst-case multi-objective error estimate can be decomposed and bounded as:
\begin{equation}
\label{eq:errbound2}
\frac{|a_{\varepsilon}(u^{h/2}-u^h,l_{\varphi}^{h/2})|}
{\Varnorm{\smash[b]{l_{\varphi}^{h/2}}}_{H^{1/2}(\partial\Omega)}}
=
\bigg|\sum_{i\in{}\II}\sigma_ia_{\varepsilon}(u^{h/2}-u^h,\psi_i)\bigg|
\leq\sum_{i\in\II}\eta_i
\end{equation}
with $\eta_i=\big|\sigma_ia_{\varepsilon}(u^{h/2}-u^h,\psi_i)\big|$. From~\EQ{qh/2} and~\EQ{lvarphih/2} it however follows that
\begin{equation}
\label{eq:altestimates}
a_{\varepsilon}(u^{h/2}-u^h,l_{\varphi}^{h/2})
=
a_{\varepsilon}(u^{h/2}-u^h+q^{h/2},l_{\varphi}^{h/2})
=
a_{\varepsilon}(l_{\varphi}^{h/2},l_{\varphi}^{h/2})
\end{equation}
The second and third expressions in~\EQ{altestimates} provide equivalent error estimates, and can alternatively serve
to construct error bounds and error indicators analogous to~\EQ{errbound2}. It is important to note, however, that although
the error estimates in~\EQ{altestimates} are equivalent, the corresponding error bounds are generally not. From~\EQ{lvarphih/2}
one can infer that the error bounds derived from the second and third expressions in~\EQ{altestimates} coincide. However,
these will generally deviate from the error bound derived from the first expression in~\EQ{altestimates}.

\begin{figure}
\begin{center}
\epsfig{figure=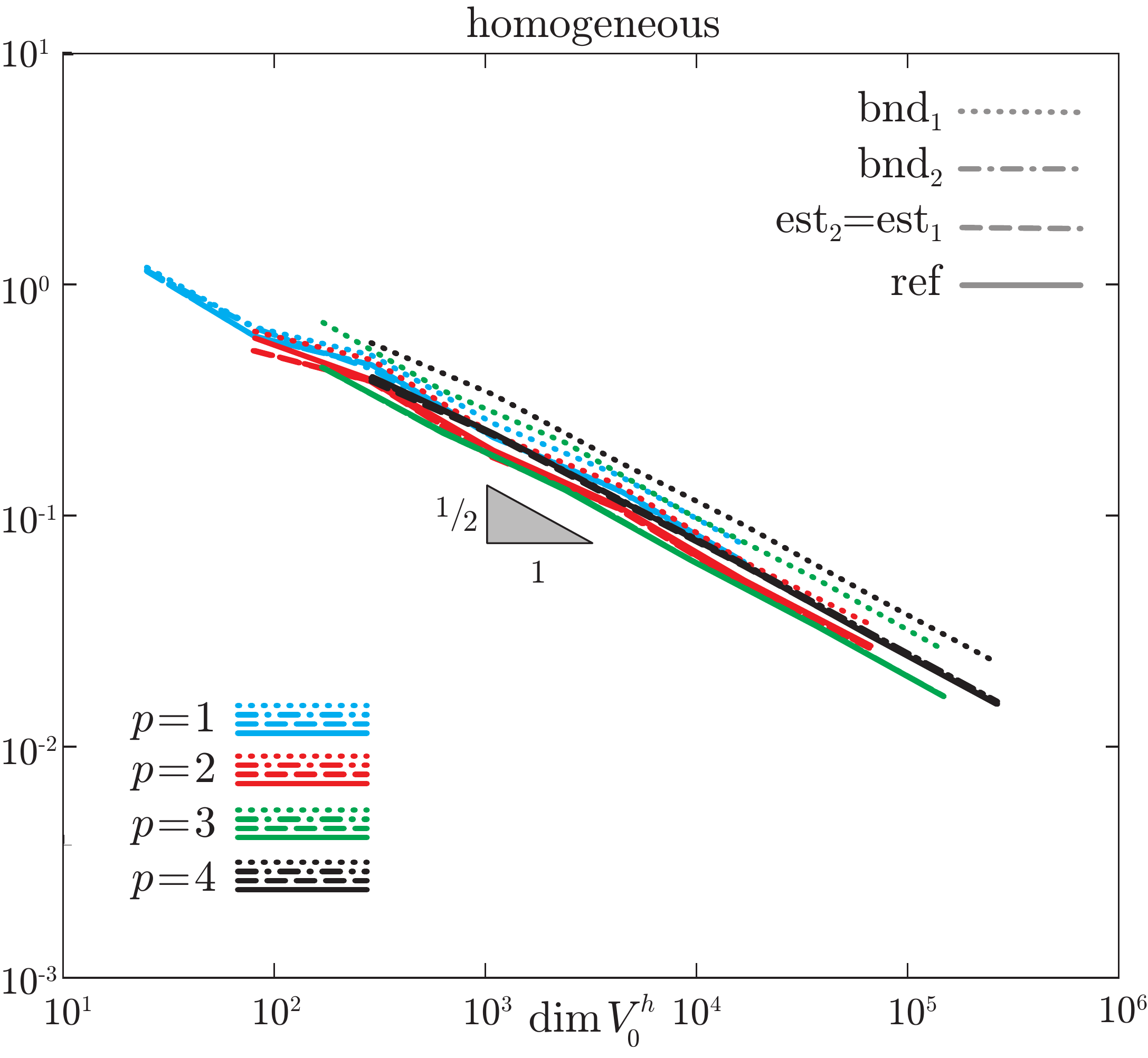,width=0.48\textwidth}
\hspace{0.02\textwidth}
\epsfig{figure=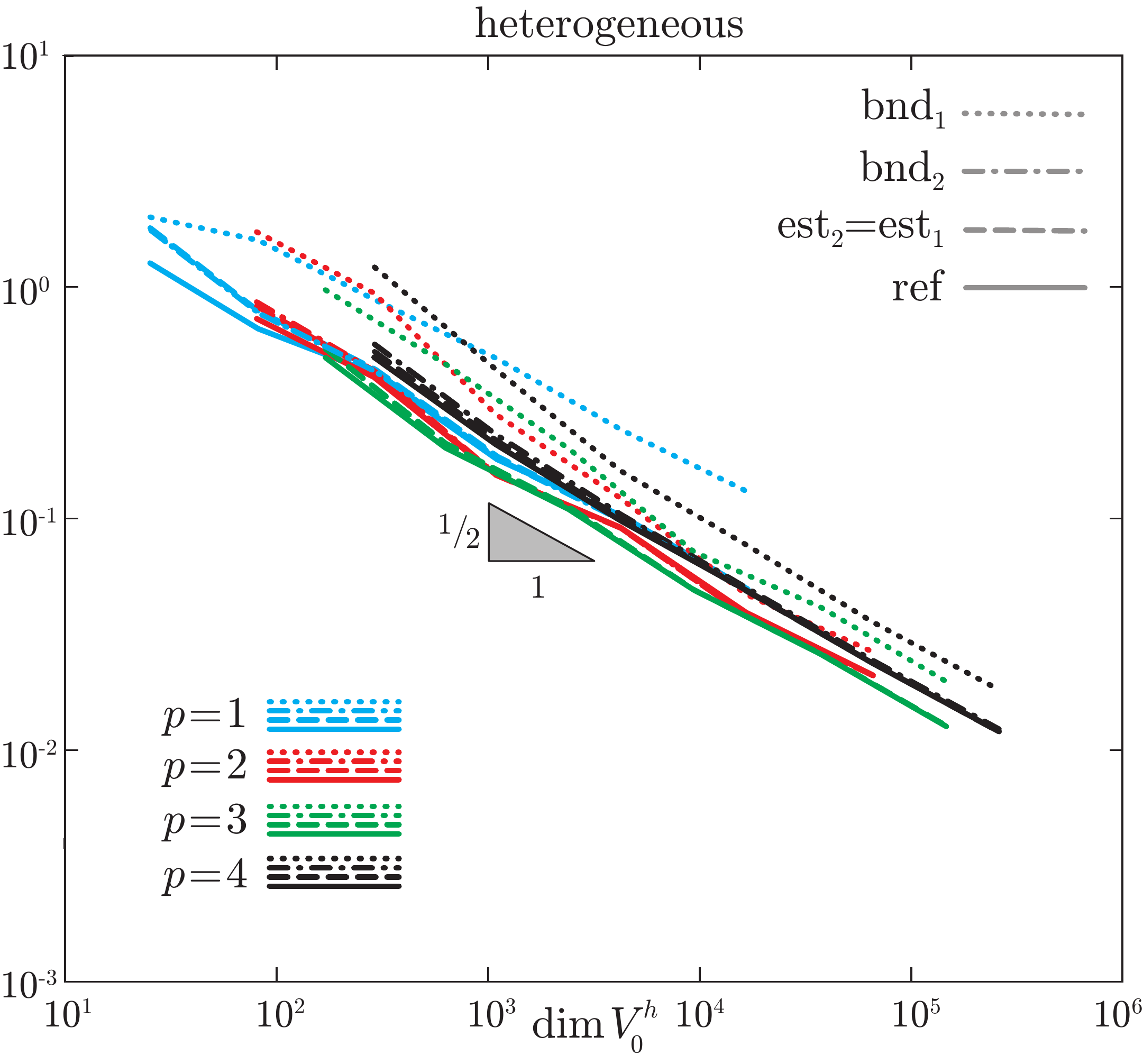,width=0.48\textwidth}
\end{center}
\caption{%
Worst-case multi-objective error estimates $\mathrm{est}_1$ and $\mathrm{est}_2$ (coincident)
for test case~3 with homogeneous coefficients ({\em left\/}) and heterogeneous coefficient ({\em right\/}),
and upper bounds
$\mathrm{bnd}_1$
and
$\mathrm{bnd}_2$
according to~\EQ{plotcontent} versus the number of degrees of freedom for finite-element spaces with polynomial orders $p\in\{1,2,3,4\}$ on uniform meshes.
The solid lines indicate the reference error estimates according to~\EQ{referror1}.
\label{fig:TC3_ConvUniform}}
\end{figure}

To assess the accuracy of the worst-case multi-objective error estimates in~\EQ{altestimates} and the tightness of the
corresponding error bounds, we consider finite-element approximations of~\EQ{weakForm} and the worst-case multi-objective
error estimates and bounds:
\begin{equation}
\label{eq:plotcontent}
\begin{aligned}
\mathrm{est}_1&=\frac{|a_{\varepsilon}(u^{h/2}-u^h,l_{\varphi}^{h/2})|}{\Varnorm{\smash[b]{l_{\varphi}^{h/2}}}_{H^{1/2}(\partial\Omega)}}
\\
\mathrm{est}_2&=\frac{|a_{\varepsilon}(l_{\varphi}^{h/2},l_{\varphi}^{h/2})|}{\Varnorm{\smash[b]{l_{\varphi}^{h/2}}}_{H^{1/2}(\partial\Omega)}}
\end{aligned}
\quad
\begin{aligned}
\phantom{\frac{l_{\varphi}^{h/2}}{l_{\varphi}^{h/2}}}\mathrm{bnd}_1&=\sum_{i\in{}\II}\big|\sigma_ia_{\varepsilon}(u^{h/2}-u^h,\psi_i)\big|
\\
\phantom{\frac{l_{\varphi}^{h/2}}{l_{\varphi}^{h/2}}}\mathrm{bnd}_2&=\sum_{i\in{}\II}\big|\sigma_ia_{\varepsilon}(l_{\varphi}^{h/2},\psi_i)\big|
\end{aligned}
\end{equation}
for $\omega=\Gamma_D$, i.e. we regard estimates of the worst-case error in the flux functional along the entire Dirichlet boundary.
Figure~\FIG{TC3_ConvUniform} plots $\mathrm{est}_{1,2}$ and $\mathrm{bnd}_{1,2}$ according to~\EQ{plotcontent}
versus the dimension of the finite-element approximation space, $\smash[tb]{\dim{}V^h_{0,\Gamma_D}}$,
for standard finite-element approximations with polynomial degrees $p\in\{1,2,3,4\}$ on a sequence of uniform meshes with mesh widths
$h\in\{2^{-2},2^{-3},\ldots,2^{-7}\}$, and for the homogeneous permeability $\varepsilon=1$ ({\em left\/}) and the heterogeneous permeability
in Figure~\FIG{TC3_Setup} ({\em right\/}).
In addition, Figure~\FIG{TC3_ConvUniform} plots a reference error estimate according to
\begin{equation}
\label{eq:referror1}
\mathrm{ref}=\frac{|a_{\varepsilon}(u^{h_{\mathrm{ref}}}-u^h,l_{\varphi}^{h_{\mathrm{ref}}})|}{\Varnorm{\smash[b]{l_{\varphi}^{h_{\mathrm{ref}}}}}_{H^{1/2}(\partial\Omega)}}
\end{equation}
with $h_{\mathrm{ref}}=2^{-8}$. Comparison of the solid and dashed curves in Figure~\FIG{TC3_ConvUniform} indicates that the error
estimates $\mathrm{err}_{1,2}$ generally provide accurate estimates of the worst-case multi-objective error, even on coarse
meshes. Only for the heterogeneous case, $p=1$ and coarse meshes ($h\in\{2^{-2},2^{-3}\}$) there is a noticeable deviation between $\mathrm{err}_{1,2}$
and the reference error estimate.
Figure~\FIG{TC3_ConvUniform} moreover reveals a significant difference between $\mathrm{bnd}_1$ and $\mathrm{bnd}_2$, despite the fact
that these bounds derive from equivalent error estimates. While $\mathrm{bnd}_2$ is tight and the deviation between the
estimates $\mathrm{err}_{1,2}$ and the bound~$\mathrm{bnd}_2$ is generally not discernible, the bound $\mathrm{bnd}_1$ typically deviates from $\mathrm{err}_{1,2}$
by a factor of 2 or more. The deviation between $\mathrm{bnd}_1$ and $\mathrm{err}_{1,2}$ is particularly manifest for the heterogeneous-permeability case.
One may also observe from Figure~\FIG{TC3_ConvUniform} that both in the homogeneous and heterogeneous case,
the error decays only as $(\smash[tb]{\dim{}V^h_{0,\Gamma_D}})^{1/2}$, independent of the order of approximation. This suboptimal convergence
behavior is attributable to the fact that the near-singularities in the Dirichlet data are not  resolved on the considered meshes.


%
\subsection{Worst-case multi-objective error estimation and adaptivity with data incompatibility}
\label{sec:numer-exper-multi-object-adapt-incompatible}
In this section we consider worst-case multi-objective error estimation and adaptivity for weak formulation~\EQ{weakForm}
in conjunction with objective set~\EQ{TC3O}. The adaptive-refinement procedure based on the worst-case multi-objective error estimate
associated with~\EQ{TC3O} essentially aims to construct a finite-element space that yields an optimal approximation of the
flux functional $-\varepsilon\partial_nu$ on~$\omega=\Gamma_D$.
We again consider the bilinear form~\EQ{aeps} with both a homogeneous coefficient $\varepsilon=1$ and a heterogeneous coefficient
according to Figure~\FIG{TC3_Setup}.
Recalling that $\mathrm{bnd}_2$ in~\EQ{plotcontent} generally provides a tighter bound of the error estimate than~$\mathrm{bnd}_1$, we base the marking strategy
in the SEMR procedure on the error indicators $\eta_i=|\sigma_ia_{\varepsilon}(l_{\varphi}^{h/2},\psi_i)|$,
i.e. the summands in $\mathrm{bnd}_2$. The marking and refinement
operations are otherwise identical to those described in Section~\SEC{numer-exper-multi-object-adapt-compatible}. The initial mesh for the adaptive computations is composed of $4\times{}4$ elements.

\begin{figure}
\begin{center}
\centering\epsfig{figure=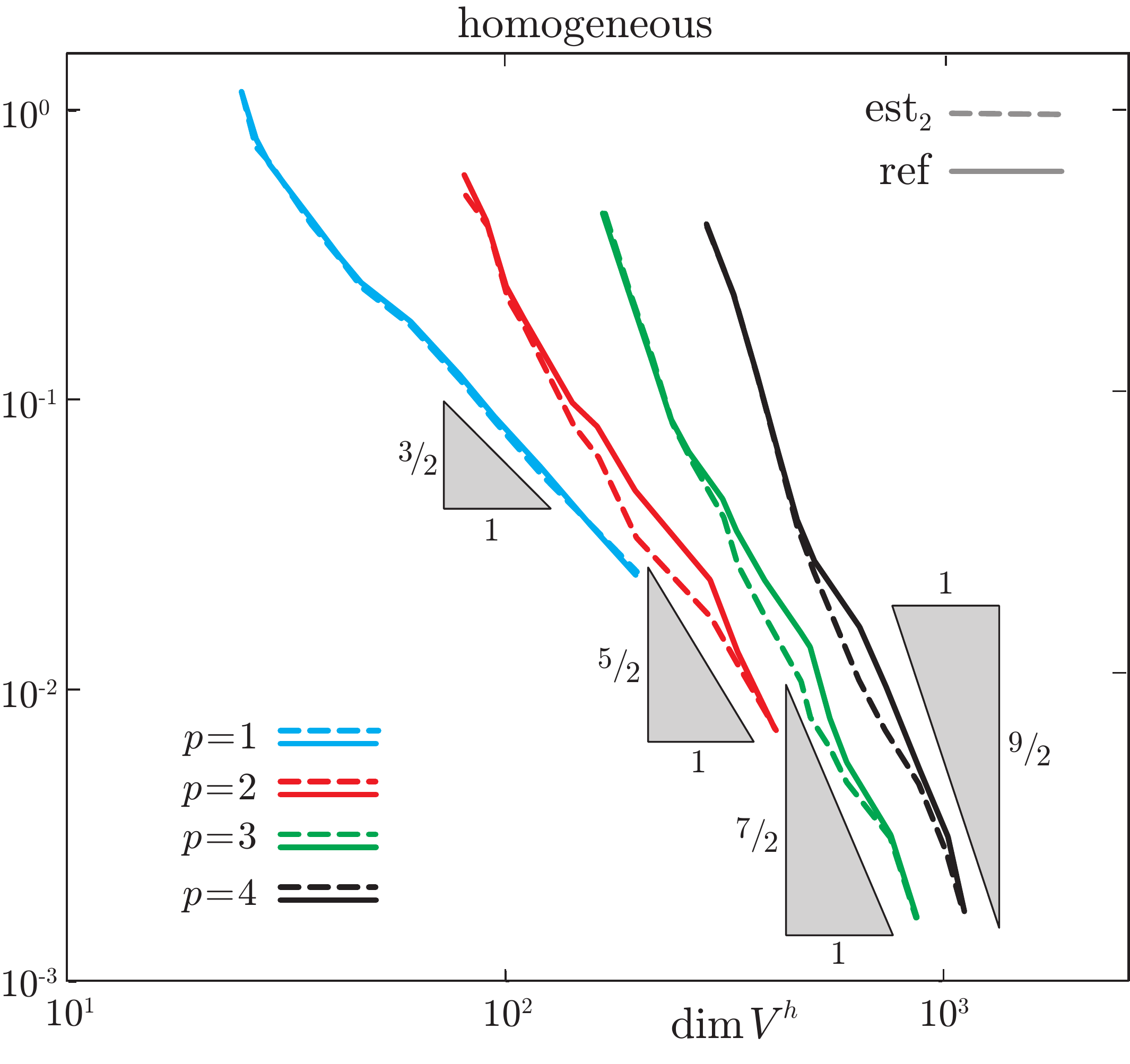,width=0.48\textwidth}
\hspace{0.02\textwidth}
\centering\epsfig{figure=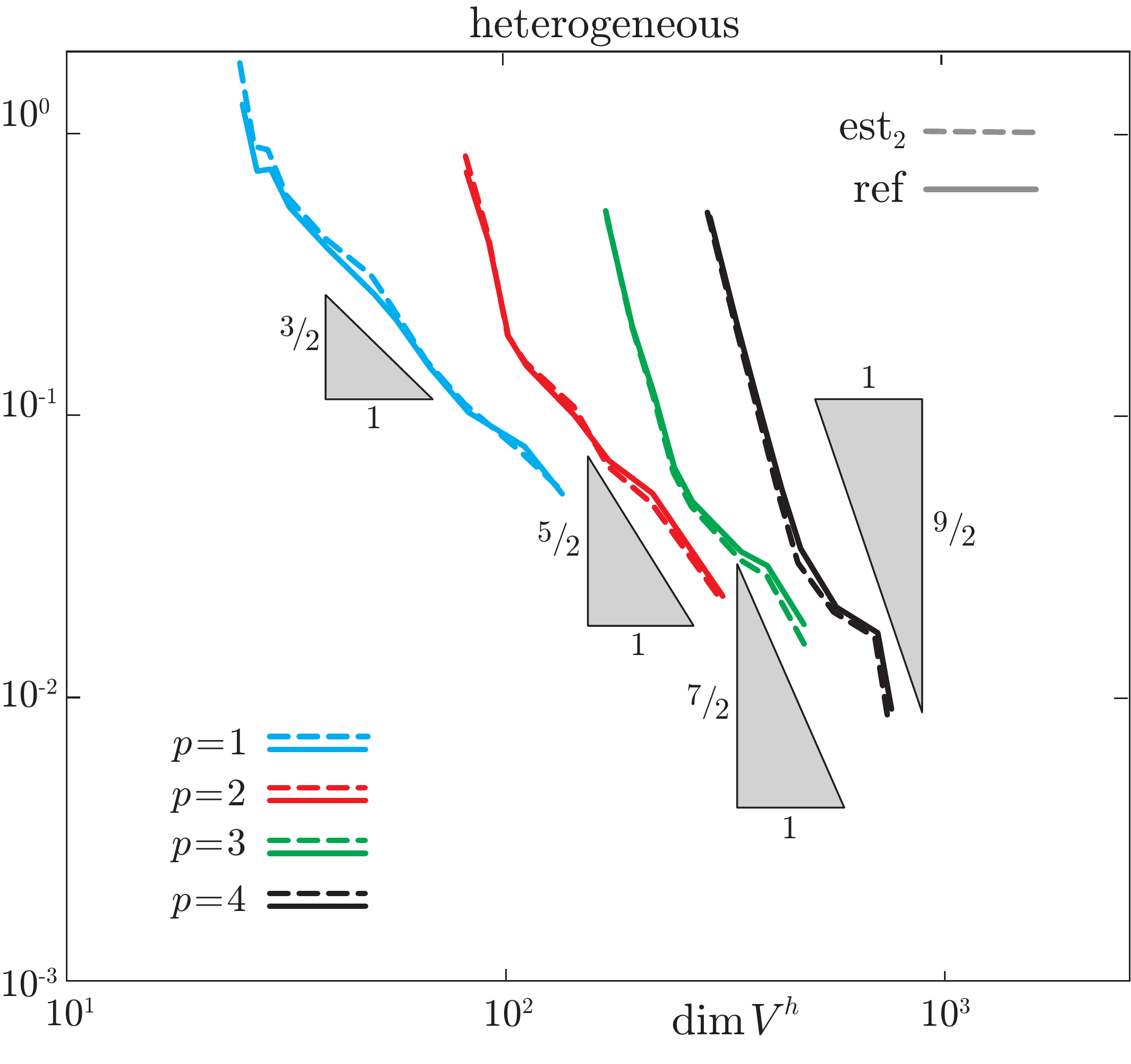,width=0.48\textwidth}
\caption{%
Worst-case multi-objective error estimate $\mathrm{est}_2$ 
and reference error estimate 
for test case~3 with homogeneous coefficients ({\em left\/}) and heterogeneous coefficient ({\em right\/}),
versus $\smash[tb]{\dim{}V^h}$ for sequences of adaptive finite-element approximations with polynomial orders $p\in\{1,2,3,4\}$.
\label{fig:TC3_ConvAdapt}}
\end{center}
\end{figure}

Figure~\FIG{TC3_ConvAdapt} displays the worst-case multi-objective error estimate $\mathrm{est}_2$ according to~\EQ{plotcontent}
corresponding to the finite-element approximation~$u^h$ generated by the SEMR algorithm, versus the dimension of the adaptively
refined finite-element space for finite-element approximations with polynomial orders $p\in\{1,2,3,4\}$.
Results for the homogeneous coefficient $\varepsilon=1$ and the heterogeneous coefficient according to Figure~\FIG{TC3_Setup} are
presented in the left and right panels of Figure~\FIG{TC3_ConvAdapt}, respectively. In addition to the error estimate $\mathrm{est}_2$,
Figure~\FIG{TC3_ConvAdapt} displays a reference error estimate conforming to~\EQ{referror1} with $V^{h_{\mathrm{ref}}}$ the finite-element
space that is obtained by uniform bisection of all elements in the final mesh constructed by the adaptive algorithm. It is noteworthy that the deviation between the error estimate and the reference
error estimate is generally small. The
left panel of Figure~\FIG{TC3_ConvAdapt} indicates that in the homogeneous case, the
adaptive-refinement procedure yields a convergence rate proportional to~$(\mathrm{dim}V^h)^{-(p+1/2)}$. A similar convergence rate
appears to hold for the heterogeneous case, although the results in the right panel of Figure~\FIG{TC3_ConvAdapt} are
less conclusive. It is to be noted, however, that for the heterogeneous case the refinement procedure reaches the resolution of the permeability data ($628\times628$) after 7 levels of refinement, which affects the asymptotic convergence rate. We conjecture that the convergence rate~$(\mathrm{dim}V^h)^{-(p+1/2)}$ corresponds
to optimal convergence in the $H^{1/2}$\nobreakdash-norm in 1~dimension, but such approximability results in fractional-order
Sobolev spaces have to our knowledge only been established for piecewise-linear approximations~\cite{Ciarlet-Jr.:2013dn}.
It is remarkable that the convergence rate achieved by the adaptive algorithm appears to
pertain to approximation of a 1\nobreakdash-dimensional object.
This is however consistent with the fact that the worst-case multi-objective error corresponds to
$s(\DD,u-u^h)=\varnorm{\varphi}_{H^{1/2}(\partial\Omega)}$ with $\smash[tb]{\varphi\in{}H^{1/2}_{0,\partial\Omega\setminus\omega}(\partial\Omega)}$; see~\EQ{TC3ids} and~\EQ{varphi1}.
It is therefore plausible that (for each $p$) the adaptive algorithm generates a
sequence of adaptive finite-element approximation spaces, $\{V^h\}$, such that
the worst-case multi-objective error of the Galerkin approximations $u^h\in{}V^h$
decays as $s(\DD,u-u^h)=\varnorm{\varphi}_{H^{1/2}(\partial\Omega)}\propto(\mathrm{dim}V^h)^{-(p+1/2)}$ as
$\mathrm{dim}V^h\to\infty$. A detailed analysis of this
convergence behavior is beyond the scope of this work.

\begin{figure}
\begin{center}
\centering\epsfig{figure=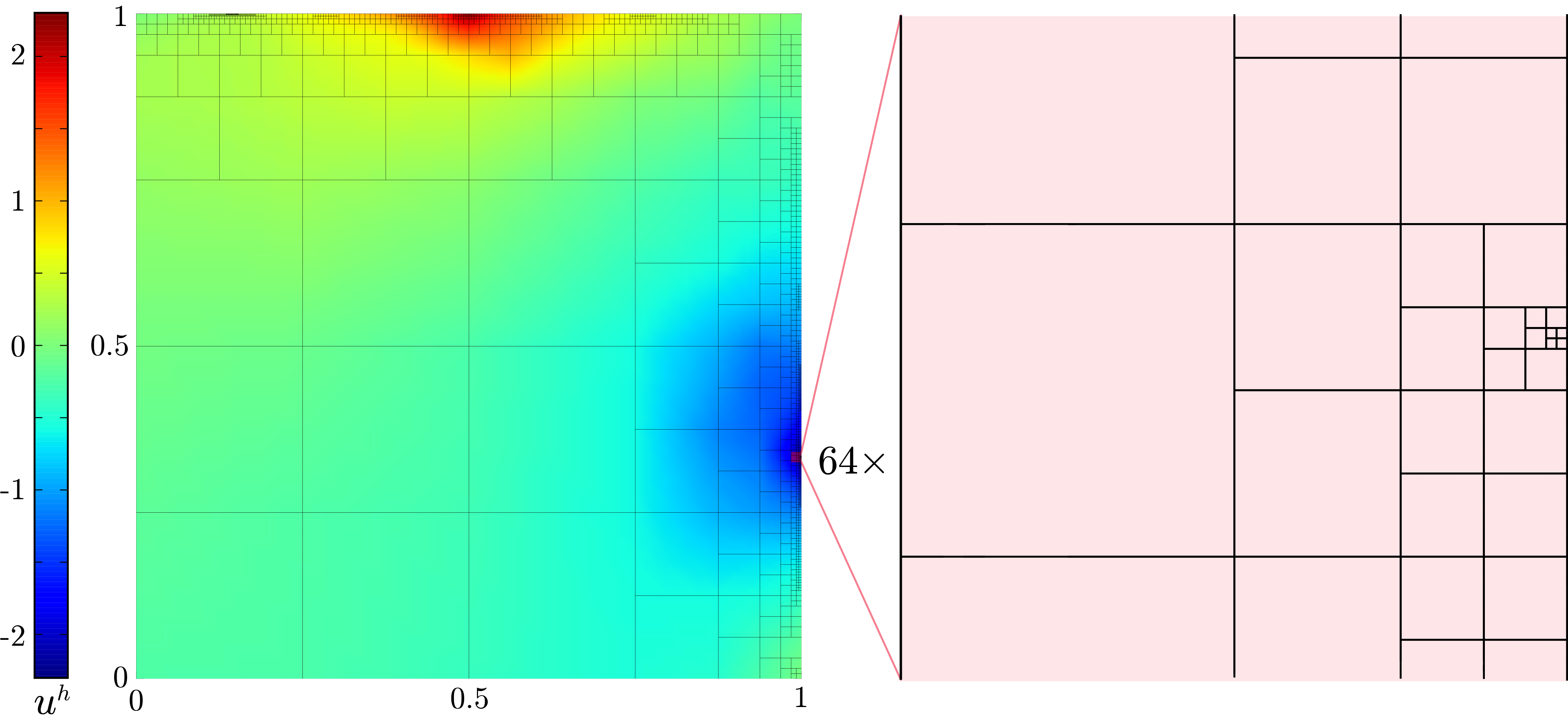,width=0.99\textwidth}
\caption{%
Finite-element mesh for test case~3 with heterogeneous coefficients for approximation order $p=1$,
after 16 iterations of the SEMR process. The red square presents a
magnification ($64\times$) of the mesh near the nearly-singular point $(x,y)=(1,1/3)$.
\label{fig:TC3_Mesh}}
\end{center}
\end{figure}

Figure~\FIG{TC3_Mesh} presents the finite-element mesh that is obtained after 16 iterations of the SEMR process for the test case with a heterogeneous coefficient, for a finite-element approximation of order $p=1$.
The red square in Figure~\FIG{TC3_Mesh}
presents a magnification ($64\times$) of the mesh near the nearly-singular point $(x,y)=(1,1/3)$.
One can observe that the mesh is refined towards the part of the boundary corresponding to $\omega=\Gamma_D\subset\partial\Omega$ which constitutes the support of the functional $-\varepsilon\partial_nu$
according to the objective set~\EQ{TC3O}. Moreover, the adaptive refinement is most concentrates near the nearly-singular points and, in particular, near the irrational point $(x,y)=(1,1/3)$. Comparing the mesh in
Figure~\FIG{TC3_Mesh} with the value of the heterogeneous coefficent $\varepsilon$ in Figure~\vref{fig:TC3_Setup}, one can notice that the adaptive mesh is generally finer in regions where the coefficient (permeability) is relatively large. Let us also allude to the refinement near the bottom
right corner, $(x,y)=(1,0)$.
It appears that a weak singularity occurs at this point due the transition between the Dirichlet
boundary condition and the Neumann boundary condition.

\section{Conclusion}
\label{sec:Concl}
In this work we introduced a new computational methodology for determining
a-posteriori worst-case multi-objective error estimates for finite-element approximations.
The methodology applies to both standard finite-element approximation errors and data-incompatibility errors
due to incompatibility of boundary data with the trace of the finite-element space.
As opposed to goal-oriented approaches, which consider only a single objective functional, the presented methodology applies to general closed convex subsets of the dual space.
The worst-case multi-objective error coincides with the support function of the considered objective set
at the finite-element approximation error. The error estimate presented in this work
adopts a standard dual-weighted-residual form, in which the dual solution corresponds
to an approximation of the supporting functional of the objective set at the approximation error.
The error estimate can direct an adaptive-refinement procedure in a similar manner as in conventional goal-oriented approaches, rendering a multi-objective adaptive strategy.

To illustrate the properties of the worst-case multi-objective error-estimation technique and its application in adaptive refinement procedures, we presented numerical results for two Dirichlet--Poisson test cases with compatible boundary data and one Dirichlet--Neumann--Laplace test case with incompatible boundary data. In the numerical experiments, we considered an approximation of the supporting functional based on the finite-element approximation and a finite-element approximation on a mesh refined by uniform bisection. We generally observed very good agreement between the worst-case multi-objective error estimate thus obtained and the actual worst-case multi-objective error, even on coarse meshes and at low orders of approximation.

For the test case with incompatible boundary data, we derived two distinct but equivalent error estimates. We showed that in contrast to the equivalence of the error estimates, the error bounds obtained from these estimates are generally not equivalent. While one error bound was observed to be tight, the other displayed a significant deviation from the underlying error estimate.

For the considered test cases, application of the worst-case multi-objective error estimate in an adaptive refinement procedure resulted in optimal convergence rates. For the Dirichlet--Neumann--Laplace test case with data incompatibility, the adaptive algorithm yields a form of super-convergence, in that
the obtained convergence rate corresponds to an optimal convergence rate for a function on a manifold of co-dimension 1.

\bibliography{BibFile}
\bibliographystyle{plain}

\end{document}

%% file: preamble.tex
\usepackage{multicol}
\usepackage{multirow}
\usepackage{fancybox}
\usepackage{index}
\usepackage{varioref}
\usepackage{psfrag}
\usepackage{epsfig}
\usepackage{boxedminipage}
\usepackage{graphicx}
\usepackage{rotating}
\usepackage{amsmath}
\usepackage{amssymb}
\usepackage{latexsym}
\usepackage{alltt}
\usepackage[small,bf]{caption}
\usepackage{url}
\usepackage{citesort}
\usepackage{array}
\usepackage{subfigure}
\usepackage{dcolumn}
\usepackage{mathrsfs}

\newcommand{\nc}{\newcommand}
\nc{\mathsm}[1]{\text{\small{$#1$}}}
\nc{\ubar}[1]{\underset{-}{#1}}
\nc{\optype}{\textrm}
\nc{\EQ}[1]{(\ref{eq:#1})}
\nc{\TAB}[1]{\ref{tab:#1}}
\nc{\FIG}[1]{\ref{fig:#1}}
\nc{\SEC}[1]{\ref{sec:#1}}
\nc{\CHAP}[1]{\ref{chap:#1}}
\nc{\THEO}[1]{\ref{theo:#1}}
\nc{\mtrx}[1]{\boldsymbol{\mathbf{#1}}}
\nc{\DD}{\mathcal{O}}
\nc{\OO}{\mathcal{O}}
\nc{\vctr}[1]{{#1}}
\nc{\grad}{\nabla}

\newcommand{\closeremark}{\hfill{\small $\blacksquare$}}
\newcommand{\jsup}{\jmath_{\mathrm{s}}}
\nc{\gradient}{\textsl{grad}\,}
\nc{\hessian}{\textsl{grad\,}^2}
\nc{\ii}{\iota}
\nc{\dd}{\mathrm{d}}
\nc{\ee}{e}
\nc{\pdiv}[2]{\partial{#1}/\partial{#2}}
\nc{\dpdiv}[2]{\displaystyle{\frac{\partial{#1}}{\partial{#2}}}}
\nc{\ddiv}[2]{\displaystyle{\frac{\dd{#1}}{\dd{#2}}}}
\nc{\inpr}{\hspace{-1pt}\cdot\hspace{-1pt}}
\nc{\IR}{\mathbb{R}}
\nc{\IN}{\mathbb{N}}
\nc{\IZ}{\mathbb{Z}}
\nc{\IC}{\mathbb{C}}
\nc{\half}{\frac{1}{2}}
\nc{\shalf}{\scriptstyle{\half}} 
\nc{\ds}[1]{\displaystyle{#1}}
\nc{\ts}[1]{\textstyle{#1}}
\nc{\sign}{\optype{sign}}
\nc{\bdy}{\optype{bdy}}
\nc{\interior}{\optype{int}}
\nc{\spr}{\optype{spr}}
\nc{\dist}{\optype{dist}}
\nc{\cond}{\optype{cond}}
\nc{\kernel}{\optype{kernel}}
\nc{\spa}{\optype{span}}
\nc{\order}{\mathcal{O}}
\nc{\Fr}{\mathrm{Fr}}
\nc{\Rey}{\mathrm{Re}}
\nc{\Ord}{O}
\nc{\ord}{o}
\nc{\st}{\:{:}\:}
\nc{\closure}[1]{\overline{#1}}
\nc{\emin}[1]{\emph{#1}\index{#1}\/}
\nc{\rmin}[1]{#1\index{{}@{#1}}}
\nc{\Laplace}{\Delta}
\nc{\ie}{i.e.}
\nc{\eg}{e.g.}
\nc{\union}{\cup}
\nc{\Union}{\bigcup}
\nc{\lf}[1]{\mathsf{#1}}
\nc{\dbar}[1]{\bar{\bar{#1}}}
\nc{\ul}[1]{\underline{#1}}
\nc{\hpt}{\hspace{0.5pt}}
\nc{\E}[1]{\times{}10^{#1}}
\nc{\inp}[2]{\langle{#1},{#2}\rangle}
\nc{\tmpcommand}{}
\nc{\GG}{{G}}
\nc{\G}[1]{\GG({#1})}
\nc{\Par}[1]{\mathscr{G}({#1})}
\nc{\const}{\mathscr{C}}
\newcommand{\varnorm}[1]{{\left\vert\kern-0.25ex\left\vert\kern-0.25ex\left\vert #1 
    \right\vert\kern-0.25ex\right\vert\kern-0.25ex\right\vert}}
\newcommand{\Varnorm}[1]{{\big\vert\kern-0.25ex\big\vert\kern-0.25ex\big\vert #1 
    \big\vert\kern-0.25ex\big\vert\kern-0.25ex\big\vert}}
    
\newcommand{\IA}{\mathcal{A}}
\newcommand{\II}{\mathcal{I}}
\newcommand{\IM}{\mathcal{M}}

\DeclareMathOperator*{\argmin}{arg\:min}

\newcommand{\meas}[1]{\optype{meas}}
\renewcommand{\half}{\mbox{$\frac{1}{2}$}}